\documentclass[12pt]{article}

\usepackage{amsmath,amssymb,amsthm,amstext,amsfonts,mathrsfs,latexsym}
\usepackage{txfonts}
\usepackage{graphicx}
\usepackage{tikz}
\usetikzlibrary{arrows.meta,calc,positioning}
\usepackage{upgreek}
\usepackage{multirow}
\usepackage{tabularx}
\usepackage{float}
\usepackage{authblk}
\usepackage{url}
\usepackage{color}
\usepackage[left=1.4in,right=1.4in,top=1.4in,bottom=1.4in]{geometry}
\newcommand{\R}{\mathbb{R}}

\DeclareMathOperator{\tr}{tr}

\newcommand{\D}{\mathcal{D}}
\newcommand{\G}{\mathcal{G}}
\newcommand{\pa}{\mathrm{pa}}

\providecommand{\Perp}{}
\renewcommand{\Perp}{\perp\!\!\!\perp}

\newcommand{\Pg}{\mathrm{P}_{\G}}
\newcommand{\Pd}{\mathrm{P}_{\D}}

\newcommand{\Rd}{\mathrm{R}_{\D}}
\newcommand{\Sd}{\mathrm{S}_{\D}}

\tikzset{
  dag vertex/.style={circle,draw=black,fill=white,minimum size=5mm,
    inner sep=0pt,font=\scriptsize},
  dag edge/.style={-{Stealth[length=1.7mm,width=1.1mm]},line width=0.5pt},
  molecular true edge/.style={-{Stealth[length=1.3mm,width=0.9mm]},
    blue!72!black,line width=0.45pt,shorten <=2.4pt,shorten >=3.1pt},
  molecular false edge/.style={-{Stealth[length=1.3mm,width=0.9mm]},
    red!78!black,densely dashed,line width=0.45pt,
    shorten <=2.4pt,shorten >=3.1pt},
  molecular vertex/.style={circle,draw=black!60,fill=cyan!32,
    minimum size=3.6mm,inner sep=0pt,line width=0.35pt},
  molecular label/.style={font=\scriptsize,text=blue!45!black,
    label distance=0.2mm}
}

\newcommand{\SachsNetworkPanel}[3]{%
  \begin{minipage}[t]{0.32\textwidth}
  \centering
  {\scriptsize\bfseries #1\par}\vspace{1mm}
  \resizebox{\linewidth}{!}{%
    \begin{tikzpicture}[x=1cm,y=1cm]
      \coordinate (cP38)  at (0,2.15);
      \coordinate (cJNK)  at (-1.18,1.80);
      \coordinate (cPLCG) at (1.18,1.80);
      \coordinate (cAKT)  at (-1.98,0.85);
      \coordinate (cPKA)  at (1.98,0.85);
      \coordinate (cERK)  at (-2.15,-0.40);
      \coordinate (cPIP3) at (2.15,-0.40);
      \coordinate (cMEK)  at (-1.65,-1.55);
      \coordinate (cPIP2) at (1.65,-1.55);
      \coordinate (cRAF)  at (-0.60,-2.23);
      \coordinate (cPKC)  at (0.60,-2.23);
      \foreach \from/\to in {#2}
        \draw[molecular true edge] (c\from) -- (c\to);
      \if\relax\detokenize{#3}\relax\else
        \foreach \from/\to in {#3}
          \draw[molecular false edge] (c\from) -- (c\to);
      \fi
      \foreach \name in {P38,JNK,PLCG,AKT,PKA,ERK,PIP3,MEK,PIP2,RAF,PKC}
        \node[molecular vertex,
          label={[molecular label]below:{\name}}] at (c\name) {};
    \end{tikzpicture}%
  }
  \end{minipage}%
}

\newcommand{\ColliderDAG}{%
  \begin{tikzpicture}[x=1cm,y=1cm]
    \coordinate (colliderthree) at (-1.05,0.70);
    \coordinate (collidertwo) at (1.05,0.70);
    \coordinate (colliderone) at (0,-0.75);
    \draw[dag edge,shorten <=2.5mm,shorten >=2.8mm]
      (colliderthree) -- (colliderone);
    \draw[dag edge,shorten <=2.5mm,shorten >=2.8mm]
      (collidertwo) -- (colliderone);
    \node[dag vertex] at (colliderthree) {$3$};
    \node[dag vertex] at (collidertwo) {$2$};
    \node[dag vertex] at (colliderone) {$1$};
  \end{tikzpicture}%
}

\theoremstyle{plain}
\newtheorem{theorem}{Theorem}[section]
\newtheorem{lemma}[theorem]{Lemma}
\newtheorem{proposition}[theorem]{Proposition}
\newtheorem{corollary}[theorem]{Corollary}

\theoremstyle{definition}
\newtheorem{definition}[theorem]{Definition}

\newtheorem{Ex}[theorem]{Example}
\newtheorem{algo}[theorem]{Algorithm}

\theoremstyle{remark}
\newtheorem{Rem}[theorem]{Remark}

\title{High dimensional Bayesian inference for Gaussian directed acyclic graph models}
\author[1]{Emanuel Ben-David \thanks{bendavid@stat.columbia.edu}}
\author[2]{Tianxi Li \thanks{tianxili@umich.edu}}
\author[3]{H\'{e}l\`{e}ne Massam\thanks{massamh@mathstat.yorku.ca}}
\author[4]{Bala Rajaratnam \thanks{brajarat@stanford.edu}}
\affil[1]{Department of Statistics, Columbia University}
\affil[2]{Department of Statistics, University of Michigan}
\affil[3]{Department of Mathematics and Statistics, York University}
\affil[4]{Department of Statistics, Stanford University}
\date{Corrected working preprint v6.2, August 2026}

\begin{document}
\maketitle

\begin{abstract}
We study centered Gaussian models Markov with respect to a directed
acyclic graph (DAG) whose vertices have a fixed parent ordering.  We
construct a conjugate family on the modified Cholesky parameters, with
one shape parameter per vertex, and derive its induced distributions on
incomplete covariance and precision coordinates.  The distribution is
proper exactly when \(\alpha_i>|\mathrm{pa}(i)|+2\) for every vertex,
and its nodewise conditional-variance and regression parameters are
independent across vertices.  This factorization gives a closed-form
normalizing constant, conjugate updating, marginal likelihoods, and
explicit full-matrix posterior means.  We distinguish these posterior
means from nonlinear completions of incomplete-coordinate means.  We also
distinguish the transformed Cholesky-coordinate mode from modes defined
using covariance or precision coordinates.  The
model-selection procedure searches only over DAGs compatible with the
specified ordering.  Historical simulation and data examples
illustrate the method; their evidentiary limitations and reproducibility
requirements are stated explicitly.
\end{abstract}

\begin{center}
\begin{minipage}{0.92\textwidth}
\small\textbf{Revision note.}
This working correction preserves the model and main construction of the
2015 arXiv version while repairing non-integrable displayed kernels,
sign inconsistencies, covariance maps, posterior-mean and mode
interpretations, and the Hausdorff-measure example.  Version 6.1 also
removes an invalid scalar density attributed to the inverse type-II
Wishart and replaces the asserted family identity by the precise
block-factorization relationship supported by \cite{letac07*}.  Empirical
claims that cannot be reproduced from the public artifacts are labeled as
historically reported rather than independently verified.  Version 6.2
replaces every externally rendered DAG figure, including all three
molecular-network panels, by native TikZ source embedded in this document.
\end{minipage}
\end{center}

%%%%%%%%%%%%%%%%
%%%%%%%%%%%%%%%%%%%%%%%%%%%%%%%%%%%%%%%  Itroducton 

\section{Introduction}\label{sec:intro}
Priors for Gaussian distributions Markov with respect to a DAG have a
long history beginning with the compatible parameter priors studied by
Geiger and Heckerman \cite{Geiger1994,Heckerman02}.  Dawid and Lauritzen
\cite{dawid93*} introduced the hyper inverse Wishart for decomposable
undirected graphs, and Letac and Massam \cite{letac07*} developed the
richer type-I and type-II Wishart families on the corresponding cones.
For a perfect orientation of a decomposable graph, Theorem~4.4 of
\cite{letac07*} gives a block inverse-Wishart/matrix-normal factorization
of the inverse type-II law.  This is the decomposable-graph analogue of
the vertexwise conditional-variance/regression factorization below.

The relationship is a specialization, not an unrestricted identity of
families.  The type-II shapes are attached to clique and separator blocks
and satisfy admissibility constraints, whereas the present construction
uses one free shape parameter per vertex.  Compatible tied choices recover
the classical Wishart and hyper-inverse-Wishart special cases, and certain
perfect-DAG specializations correspond to inverse type-II laws.  Arbitrary
vertexwise DAG--Wishart shapes need not admit a type-II parameterization.
Our construction therefore starts directly on modified-Cholesky
coordinates for an arbitrary DAG and derives the covariance and precision
images afterward.  The resulting nodewise factorization gives closed-form
normalizing constants and posterior moments.

For a non-perfect DAG, the covariance and precision parameter spaces are
typically curved submanifolds of the ambient symmetric-matrix space and
therefore have ambient Lebesgue measure zero.  A distribution on either
space is nevertheless well defined as a pushforward measure; what fails
is the existence of an ambient-Lebesgue density.  We use projections
onto Euclidean spaces of functionally independent covariance and
precision entries, following the completion results in
\cite{benraj12}.  This construction yields tractable coordinate
densities.  A completion of a coordinatewise expectation must not,
however, be confused with the expectation of the completed random
matrix; Sections~\ref{sec:wishart_incomplete} and~\ref{sec:Sd} give the
full posterior-mean formulas.

The remainder is organized as follows.  Section~\ref{sec:pre} fixes the
Gaussian-DAG and modified-Cholesky notation.  Section~\ref{sec: wdag}
constructs the DAG--Wishart on $\Theta_{\D}$.
Sections~\ref{sec:wishart_incomplete} and~\ref{sec:Sd} derive its images on
incomplete precision and covariance coordinates, including the corrected
full-matrix moments.  Section~\ref{sec:empirical} preserves the fixed-order
simulation and data examples while separating reproducible facts from
historical numerical reports.  Section~\ref{sec:closing} states the main
limitations; Supplemental Sections A--D provide proofs, geometry, and
algorithms.
%%%%%%%%%%%%%%%%%%%%                            Preliminaries
%%%%%%%%%%%%%%%%
%%%%%%%%%%%%%%%%%%%
%%%%%%%%%%%%%%%%%%%%%%%%%%%

\section{Preliminaries}\label{sec:pre}
%%%%%%%%%%%%%%%%%%%%%%%%%%%%%%%%%%%%%%%%%%%%%%%%%%%%%%%%%%%%%%%%%%%%%%

A brief summary of graph theory, associated Markov and other properties required for analyzing DAG models is given in Supplemental section A. 
% %
% %
\subsection{Gaussian DAG models}\label{sub:gdag}

Let $V$ be a set with $p$ elements. For any $a, b \subseteq V$ \footnote{Lowercase letters denote subsets of $V$.} let $\R^{a\times b}$ denote the real linear space of functions $A\equiv\left((i,j)\mapsto A_{ij}\right): a\times b\rightarrow \R$. Each element of $\R^{a\times b}$ is called an $|a|\times |b|$ matrix.  In particular, we define the space of symmetric matrices ${S}_a(\R)=\left\{A\in \R^{a\times a}: A_{ij}=A_{ji}~, \; \text{for every $i,j\in a$}\right\}$, and the set of positive definite matrices $
\mathrm{PD}_a(\R)=\left\{A\in {S}_a(\R):  ~x^{\top}Ax > 0, \; \text{for every} \: \; x\in\R^a\setminus\left\{0\right\}\right\}$.
 Now let $\Sigma\succ 0$ denote $\Sigma\in \mathrm{PD}_{p}(\R)$.  For a partition $a, b$ of $V$, consider the corresponding block partitioning of $\Sigma$ as follows.
 \[
 \Sigma=
 \left(
 \begin{matrix}
 \Sigma_{a}&\Sigma_{ab}\\
 \Sigma_{ba}&\Sigma_{b}
 \end{matrix}
 \right),
 \]
 where  $\Sigma_{aa}=(\Sigma_{ij})_{i,j\in a}\in \mathrm{PD}_a(\R)$, $\Sigma_{bb}=(\Sigma_{ij})_{i,j\in b}\in \mathrm{PD}_b(\R)$, $\Sigma_{ab}=(\Sigma_{ij})_{i\in a, j\in b}\in \R^{a\times b}$  and  $\Sigma_{ba}=\Sigma_{ab}^{\top}$. The Schur complement of the sub-matrix $\Sigma_{aa}$  is defined as  $\Sigma_{bb| a}=\Sigma_{bb}-\Sigma_{ba}(\Sigma_{aa})^{-1}\Sigma_{ab}$.
 \begin{Rem}
 Throughout this paper, we shall in general suppress the notation for a principal submatrix  $\Sigma_{aa}$ and refer to it as $\Sigma_a$. We shall also use the convention $\Sigma_a^{-1}$  for $(\Sigma_{aa})^{-1}$  and  $\Sigma^a$  for $(\Sigma^{-1})_{aa}$.
 \end{Rem}

\noindent In this paper we focus on multivariate Gaussian distributions which obey the directed Markov property with respect to  a DAG~$\D$. From now on and unless otherwise stated, we shall always assume without loss of generality that $\D=(V, E)$ is a DAG given in a parent ordering\footnote{We emphasize here that unlike in the decomposable precision graph setting or the covariance graph setting (where the existence of  an ordering is important either for the perfect order of cliques and separators, or to preserve zeros), existence of such an ordering is not necessary in the DAG setting, since a parent ordering is always available for a DAG.}, i.e., the vertices are labeled $1,2,\ldots,p$, and $i\rightarrow j$ implies that $i>j$. A Gaussian DAG model (or Gaussian Bayesian network) over $\D$, denoted by   $\mathscr{N}(\D)$,  is  the statistical model  that consists of all multivariate Gaussian distributions  $\rm{N}_{p}(\mu,\Sigma)$  obeying the ordered directed Markov property
with respect to   $\D$. Therefore, $
\mathbf{x}\sim \rm{N}_{p}(\mu,\Sigma)\in \mathscr{N}(\D)\implies x_i\Perp \mathbf{x}_{\left\{i+1,\ldots,p\right\}\setminus \pa(i)}|\mathbf{x}_{\pa(i)}$ for each $i$.
% %
 \begin{Rem} Note that $\rm{N}_{p}(\mu,  \Sigma)\in \mathscr{N}(\D)$   if and only if     $\rm{N}_{p}(0,  \Sigma)\in \mathscr{N}(\D)$. Therefore, without loss of generality, we shall only consider centered Gaussian distributions $$\left\{\rm{N}_{p}(0, \Sigma):~\Sigma\in \mathrm{PD}_{\D} \right\}\subseteq\mathscr{N}(\D).$$
\noindent For convenience, with a slight abuse of notation, we shall still denote
\[
\mathscr{N}(\D)=\left\{\rm{N}_{p}(0, \Sigma):~\Sigma\in \mathrm{PD}_{\D} \right\}.
\]
 \end{Rem}
\noindent The Gaussian distributions in $\mathscr{N}(\D)$ are naturally parametrized by the elements of 
\[
\mathrm{PD}_{\D}=\left\{\Sigma \succ 0: \: \rm{N}_{p}(0,\Sigma)\in \mathscr{N}(\D)\right\}\text{ or } \:  \mathrm{P}_{\D}=\left\{\Sigma^{-1}\succ 0: \:   \rm{N}_{p}(0,\Sigma)\in \mathscr{N}(\D)\right\}.
\]
 These sets are referred to as the space of covariance matrices and the space of precision matrices. A precision matrix in $\Pd$ is usually denoted by $\Omega$. Similarly, for an undirected graph $\G=(V, \mathscr{E})$ we define $\mathscr{N}(\G)$ as the set of multivariate Gaussian distributions obeying the (undirected) Markov property with respect to   $\G$. In this model the corresponding parameter spaces are the space of covariance matrices $ \mathrm{PD}_{\G}=\left\{\Sigma:  \rm{N}_{p}(0,\Sigma)\in \mathscr{N}(\G)\right\}$ and the space of precision matrices $\mathrm{P}_{\G}=\left\{\Omega :  \Omega^{-1}\in \mathrm{PD}_{\G}\right\}$. Note that, for us, $\mathrm{PD}_{\D}$ and $\Pd$ are  parameter spaces of primary interest as they arise naturally in the parameterization of  Gaussian densities. However, in order to develop multi-shape parameters Wishart priors on these spaces, which is the main purpose of this paper,we begin with the more natural and more convenient  Cholesky type parameterization of $\mathscr{N}(\D)$ that we discuss in the next subsection.  
\subsection{Cholesky parametrizations of Gaussian DAG models}\label{subs:d_parameterization}
Consider a Gaussian DAG distribution $\rm{N}_{p}(0,\Sigma )\in
\mathscr{N}(\D)$.  Let $\mathcal{L}_{\D}$ denote the set of lower
triangular matrices with unit diagonal and \(L_{ij}=0\) if
\(i\notin\mathrm{pa}(j)\), and let $\mathcal{D}_{+}^{p}$ denote the set
of strictly positive diagonal matrices in $\R^{p\times p}$.  Then
\(\Sigma^{-1}\in\Pd\) if and only if there is a unique
\((D,L)\in\mathcal{D}_{+}^{p}\times\mathcal{L}_{\D}\) such that
\[
\Omega=\Sigma^{-1}=LD^{-1}L^{\top},
\qquad
\Sigma=L^{-\top}DL^{-1}.
\]
The first identity is the modified Cholesky decomposition of
\(\Omega\).  We call
\(\Theta_{\D}=\mathcal{D}_{+}^{p}\times\mathcal{L}_{\D}\) the Cholesky
space of \(\D\), and
\[
\left\{\mathrm{N}_p(0,L^{-\top}DL^{-1}): (D,L)\in\Theta_{\D}\right\}
\equiv\mathscr{N}(\D)
\]
the Cholesky parametrization of \(\mathscr{N}(\D)\).

\noindent We can also obtain a variant of this parameterization, in vector form, from the recursive factorization property of the Gaussian densities in  $\mathscr{N}(\D)$ (see Supplemental section A subsection 1.3 for details). First, let us recall the following notation from \cite{andersson98}.\\
\noindent
  \textbf{Notation.} For each $i\in V$ let
  \[
  \begin{array}{lcl}
  \prec i\succ=pa(i)&\quad[i\succ=\left\{i\right\}\times pa(i)&\quad  \prec i ]=pa(i)\times\left\{i\right\},\\
  \nprec i \nsucc=\left\{j: j>i\right\}\setminus pa(i)& \quad[ i\nsucc =\left\{i\right\}\times \nprec i \nsucc &\quad \prec i\nsucc=\prec i\succ \times \nprec i\nsucc\\
   \preceq i\succeq=fa(i) &  &
  \end{array}
  \]
  
	By applying the directed factorization property (DF) of  $\rm{N}_{p}(0, \Sigma)\in \mathscr{N}(\D)$  we  have
  \begin{align}\label{dmpogd}
	\notag d\rm{N}_{p}(0, \Sigma)(x)&=\prod_{i\in V}d\rm{N}(\mu_{i|pa(i)}, \Sigma_{ii|pa(i)})(x_{i}| x_{pa_{i}})\\
  &=\prod_{i\in V}d\rm{N}(\Sigma_{[i\succ}\Sigma^{-1}_{\prec i\succ}x_{\prec i \succ}, \Sigma_{ii|\prec i\succ})(x_i),
  \end{align}
  for each  $x=(x_i)_{i\in V}\in \R^{p}$. Note that  $\rm{N}(\Sigma_{[i\succ}\Sigma^{-1}_{\prec i\succ}x_{\prec i\succ}, \Sigma_{ii|\prec i\succ})(x_i)$  is the conditional distribution of $\mathbf{x}_{i}|\mathbf{x}_{\prec i\succ}=x_{\prec i\succ}$. Moreover, $\Sigma_{[i\succ}\Sigma^{-1}_{\prec i\succ}$ is the regression coefficient of  $\mathbf{x}_i$  in the regression of~$\mathbf{x}_i$  on  $\mathbf{x}_{\prec i\succ}$, and $\Sigma_{ii|\prec i\succ}$ is the conditional variance of $\mathbf{x}_i | \mathbf{x}_{\prec i\succ}=x_{\prec i\succ}$. Furthermore, using the exact functional form of the densities of the Gaussian distributions in \eqref{dmpogd}, we obtain the following equation.
  \begin{equation}\label{eq: Rgp}
  \tr(\Sigma^{-1}xx^{\top})=\sum_{i\in V}\tr\left(\Sigma^{-1}_{ii|\prec i\succ}(x_i-\Sigma_{[i\succ}\Sigma^{-1}_{\prec i\succ}x_{\prec i\succ}) (x_i-\Sigma_{[i\succ}\Sigma^{-1}_{\prec i\succ}x_{\prec i\succ})^{\top}\right)
  \end{equation}
  It is  shown  in \cite{andersson98} that $\Sigma\in \mathrm{PD}_{\D}$   if and only if   $\Sigma\succ 0$  and satisfies  \eqref{eq: Rgp} for all $x\in \R^{p}$. On the other hand, by the parent ordered Markov property of  $\mathbf{x}\sim \rm{N}_{p}(0, \Sigma)$ we have $\Sigma\in \mathrm{PD}_{\D}$ if $i\Perp \nprec i\nsucc |\prec  i\succ $ (or equivalently $ i\Perp \left\{i+1,\ldots, p\right\}\setminus pa(i)|pa(i)$). Hence another characterization given by \cite{andersson98} for $\Sigma\in \mathrm{PD}_{\D}$ is that $\Sigma\succ 0$  and
  \begin{equation}\label{eq:ndmp}
  \Sigma_{[i \nsucc}=\Sigma_{[i \succ}\Sigma^{-1}_{\prec  i \succ}\Sigma_{\prec i \nsucc},~ \text{for every $i\in V$.}
  \end{equation}
  We use two related vectors and do not interchange their signs:
  \[
  \ell_i:=L_{\prec i]},
  \qquad
  \beta_i:=-\ell_i
  =\Sigma_{\prec i\succ}^{-1}\Sigma_{\prec i]}.
  \]
  Thus \(\ell_i\) is the free Cholesky column and \(\beta_i\) is the
  regression coefficient of \(X_i\) on \(X_{\mathrm{pa}(i)}\).
  Defining
  $\Xi_{\D}=\times_{i\in V}\left(\R_{+}\times\R^{\prec i]}\right)$, it can be shown that the mapping
  \begin{equation} \label{eq:d_map}
  \Pi_{\D}\equiv\left(\Sigma\mapsto\times_{i\in V}( \Sigma_{ii|\prec i\succ}, \Sigma^{-1}_{\prec i\succ}\Sigma_{\prec i]})\right):  \mathrm{PD}_{\D}\rightarrow\Xi_{\D}
  \end{equation}
   is a diffeomorphism. In order to construct the inverse of this mapping let $\times_{i\in V} (\lambda_i, \beta_{\prec i]})$ denote a typical element in  $\Xi_{\D}$, with the convention that $\beta_{\prec i]}=0$ whenever $\prec i\succ = pa(i) = \emptyset$. Using \eqref{eq:ndmp}, the corresponding  $\Sigma$ can be recursively constructed starting from the largest index $p,$ by setting
\begin{equation}\label{completion_alg_QG}
\begin{cases}
  i)& \Sigma_{ii}=\lambda_i+ \beta_{\prec i]}^{\top}\Sigma_{\prec i\succ}\beta_{\prec i]};\\
	ii)& \Sigma_{\prec i]}=\Sigma_{\prec i\succ}\beta_{\prec i]};\\
	iii)&\Sigma_{[i \nsucc}=\Sigma_{[i \succ}\Sigma^{-1}_{\prec  i \succ}\Sigma_{\prec i \nsucc}.
	\end{cases}
\end{equation}
\noindent The reader is referred to \cite{andersson98} for greater detail, where in addition, it is shown that the inverse mapping above yields a positive
definite matrix in $\mathrm{PD}_{p}(\R)$, and consequently in $\mathrm{PD}_{\D}$. The mapping $\Pi_{\D}$ in \eqref{eq:d_map} gives another parametrization of $\mathscr{N}(\D)$ in terms of the elements $\times_{i\in V} \left(\lambda_i,\beta_{\prec i]} 
\right)\in \Xi_{\D}$. One can show that for each $i\in V$, $L_{\prec i]}=-\Sigma_{\prec i\succ}^{-1}\Sigma_{\prec i]}$ and $D_{ii}=\Sigma_{ii|\prec i\succ}$, therefore  each  $\times_{i\in V} \left(\lambda_i,\beta_{\prec i]} 
\right)\in \Xi_{\D}$ is, essentially, a vectorized form of a $(D,L)\in \Theta_{\D}$.
% % % % % % % % % % % % % % % % % % % % % % % % % % % % % % % % % % % % % %
% % %
% % % %
% % % % %
% % % % % % %
% % % % % % % % %

\section{The DAG-Wishart distribution on $\Theta_{\D}$}{\label{sec: wdag}} 
This section introduces a multi-shape family on the Cholesky space
$\Theta_{\D}$ by generalizing the modified-Cholesky image of an ordinary
Wishart law.  The family is defined for every DAG.  Its relationship to
the inverse type-II Wishart is through compatible perfect-DAG
specializations and the parallel directed-hyper-Markov factorization, not
through equality for arbitrary vertexwise shapes.
%
%
%
%
%%%%%%%%%%%%%%%%%%%%%%%%%%%%%%%%%%%%%%%
\subsection{DAG-Wishart densities}\label{sub:pit}
Let us start with a natural course that will lead us to the general form of the multi-shape parameter DAG-Wishart distributions on the Cholesky space $\Theta_{\D}$ with the desired properties. We begin with the classical Wishart distribution. Let us  consider $\Sigma^{-1}\sim\rm{W}_{p}(\eta, U)$ as a prior for the precision parameter of the full Gaussian model $\{ \rm{N}_{p}(0, \Sigma): \: \Sigma\succ 0\}$. Note that this model corresponds to the saturated Gaussian DAG model $\mathscr{N}(\D)$, i.e., when $\D$ is a complete DAG with $p$ vertices (see Figure \ref{fig:cdag}). Consider the mapping $\Sigma^{-1}\mapsto (D, L)$, where $(D, L)$ is the Cholesky factorization of $\Sigma^{-1}$. This mapping transforms the Wishart distribution  $\rm{W}_{p}(\eta, U)$ to a distribution on $\Theta_{\D}$ with density  proportional to
 
\begin{equation}\label{eq:thetad_c}
 \exp\left\{-\frac12\mathrm{tr}\left( \left(LD^{-1}L^{\top}\right)U\right)\right\}\prod_{i=1}^{p}D_{ii}^{-\frac{\alpha_i}{2}},
\end{equation}
with $\alpha_i=\eta+p-2i+3$.  Indeed, for the complete parent ordering,
the absolute Jacobian of \((D,L)\mapsto LD^{-1}L^\top\) is
\(\prod_{i=1}^{p}D_{ii}^{-(p-i+2)}\).  Although the
\(\alpha_i\)'s in \eqref{eq:thetad_c} look like multiple shape
parameters, they are functions of the single Wishart shape \(\eta\).
We now free these shape parameters and use the same kernel for an
arbitrary DAG.

\begin{figure}[htbp]
  \centering
  \begin{tikzpicture}[x=1cm,y=1cm]
    \coordinate (complete8) at (0,2.45);
    \coordinate (complete7) at (-1.55,1.65);
    \coordinate (complete6) at (1.55,1.65);
    \coordinate (complete5) at (-2.05,0.30);
    \coordinate (complete4) at (2.05,0.30);
    \coordinate (complete3) at (-1.50,-1.15);
    \coordinate (complete2) at (1.50,-1.15);
    \coordinate (complete1) at (0,-2.00);
    \foreach \from/\to in {
      8/7,8/6,8/5,8/4,8/3,8/2,8/1,
      7/6,7/5,7/4,7/3,7/2,7/1,
      6/5,6/4,6/3,6/2,6/1,
      5/4,5/3,5/2,5/1,
      4/3,4/2,4/1,
      3/2,3/1,
      2/1}
      \draw[dag edge,shorten <=2.5mm,shorten >=2.8mm]
        (complete\from) -- (complete\to);
    \node[dag vertex] at (complete8) {$p$};
    \node[dag vertex] at (complete7) {$p-1$};
    \node[dag vertex] at (complete6) {$p-2$};
    \node[dag vertex] at (complete5) {$p-3$};
    \node[dag vertex] at (complete4) {$p-4$};
    \node[dag vertex] at (complete3) {$3$};
    \node[dag vertex] at (complete2) {$2$};
    \node[dag vertex] at (complete1) {$1$};
  \end{tikzpicture}
  \caption{Schematic of a complete DAG in the parent ordering.  The eight
  displayed vertices illustrate the saturated orientation; intermediate
  labels are suppressed when $p>8$.}
  \label{fig:cdag}
\end{figure}
 
 \vspace{2mm}
 
For context, let $\G$ be decomposable and let $\D$ be a perfect
orientation compatible with a perfect clique order.  Theorem~4.4 of
\cite{letac07*} states that the inverse type-II Wishart is strong
directed hyper Markov: its residual covariance blocks have
inverse-Wishart laws, the associated regression blocks have conditional
matrix-normal laws, and the listed block pairs are mutually independent.
The theorem is intrinsically matrix-valued.  It does not yield the scalar
quadratic-form density printed in earlier drafts of this manuscript, and
we do not use it to prove any result for a non-perfect DAG.  The
vertexwise analogue follows directly below by completing the square in
the Cholesky density.

Let \(pa_i=|\mathrm{pa}(i)|\) and
\[
a_i=\frac{\alpha_i-pa_i}{2}-1,
\qquad
U_{ii\mid\prec i\succ}
=U_{ii}-U_{[i\succ}U_{\prec i\succ}^{-1}U_{\prec i]}.
\]
For an empty parent set, determinants below equal one and the associated
vectors and quadratic forms are absent.

\begin{definition}\label{def:dagwishart}
For \(U\in\mathrm{PD}_p(\R)\) and
\(\alpha_i>pa_i+2\) for every \(i\), the DAG--Wishart distribution on
\(\Theta_{\D}\) has density
\[
\pi_{U,\alpha}^{\Theta_{\D}}(D,L)
=\frac{1}{z_{\D}(U,\alpha)}
\exp\!\left\{-\frac12\tr(LD^{-1}L^\top U)\right\}
\prod_{i=1}^pD_{ii}^{-\alpha_i/2}
\]
with respect to \(\prod_i dD_{ii}\prod_i\prod_{j\in
\mathrm{pa}(i)}dL_{ji}\), where
\begin{equation}\label{eq:zd}
z_{\D}(U,\alpha)
=\prod_{i=1}^{p}
2^{\alpha_i/2-1}\pi^{pa_i/2}\Gamma(a_i)
\det(U_{\prec i\succ})^{-1/2}
U_{ii\mid\prec i\succ}^{-a_i}.
\end{equation}
Equivalently, using
\(\det(U_{\preceq i\succeq})
=\det(U_{\prec i\succ})U_{ii\mid\prec i\succ}\),
this is the determinant-ratio expression used in earlier versions.
\end{definition}

\begin{proposition}\label{prop:nodewise}
The kernel in Definition~\ref{def:dagwishart} is integrable if and only
if \(\alpha_i>pa_i+2\) for all \(i\).  Under the normalized law, the
pairs \((D_{ii},\beta_i)\), \(i=1,\ldots,p\), are mutually independent
and
\begin{align}
\label{eq:dii1}
D_{ii}&\sim
IG\!\left(a_i,\frac12U_{ii\mid\prec i\succ}\right),\\
\label{eq:dii2}
\beta_i\mid D_{ii}&\sim
\mathrm{N}_{pa_i}\!\left(
U_{\prec i\succ}^{-1}U_{\prec i]},
D_{ii}U_{\prec i\succ}^{-1}\right).
\end{align}
Here \(IG(a,b)\) has density
\(b^a\Gamma(a)^{-1}x^{-a-1}\exp(-b/x)\), \(x>0\).
Equivalently,
\(\ell_i\mid D_{ii}\sim
\mathrm{N}_{pa_i}(-U_{\prec i\succ}^{-1}U_{\prec i]},
D_{ii}U_{\prec i\succ}^{-1})\).
\end{proposition}

\begin{proof}
Since \(L_{\prec i]}=-\beta_i\), completing the square gives
\[
\tr(LD^{-1}L^\top U)
=\sum_{i=1}^pD_{ii}^{-1}
\left\{U_{ii\mid\prec i\succ}
+(\beta_i-U_{\prec i\succ}^{-1}U_{\prec i]})^\top
U_{\prec i\succ}
(\beta_i-U_{\prec i\succ}^{-1}U_{\prec i]})\right\}.
\]
The density therefore factorizes by vertex.  Gaussian integration over
\(\beta_i\), followed by inverse-gamma integration over \(D_{ii}\),
gives \eqref{eq:zd}; the latter integral is finite exactly when
\(a_i>0\).  The two conditional laws and mutual independence follow
from the same factorization.
\end{proof}

This nodewise factorization is the directed hyper-Markov property used
below.

Conjugacy follows from the Gaussian likelihood.  If
\[
Y_r\stackrel{\mathrm{iid}}{\sim}
\mathrm{N}_p(0,L^{-\top}DL^{-1}),\qquad
T=\sum_{r=1}^nY_rY_r^\top=nS,
\]
then
\[
(D,L)\mid Y_{1:n}\sim
\pi_{U+T,\,\alpha+n\mathbf{1}}^{\Theta_{\D}}.
\]

\begin{Rem}\label{rem:identefibility}
For fixed \(\alpha\), allowing every \(U\succ0\) is generally
nonidentifiable: the density depends only on the family blocks used in
the nodewise factorization.  Restricting \(U\) to \(\mathrm{PD}_{\D}\),
or equivalently parameterizing it by its free entries \(U^E\), restores
identifiability.  This restricted family is a full regular exponential
family when \(\D\) is perfect and has a curved representation in the
larger moral-graph natural-parameter space otherwise; Supplemental
Section~B gives the precise statement.
\end{Rem}
 % % % % % % % % % % % % % % % % % % % % % % % % %
 % % % % % %
 % % % % % % %
 % % % % % % % % %
 % % % % % % % % % % %
 % % % % % % % % % % % % % % % % % % % % % % % % % % %
 \section{The DAG-Wishart distribution on the space of incomplete precision matrices} \label{sec:wishart_incomplete}
The preceding section defined the DAG--Wishart
$\pi^{\Theta_{\D}}_{U,\alpha}$ in modified-Cholesky coordinates.  We now
derive its pushforward to incomplete precision coordinates, extending the
DAG analogue of the type-II-Wishart viewpoint beyond decomposable graphs
\cite{letac07*}.
 % % % % % % % % % % % % % % % % % % %
 % % % % % % % % % % % % % % % % % % % % % %
 % % % % % % % % % % % % % % % % % % % % % % % %
 
 \subsection{Motivation and notation}\label{subs:motive}
To follow in the tradition of Wishart-type conjugate priors, we derive
the precision and covariance images of
$\pi^{\Theta_{\D}}_{U,\alpha}$ for
$\mathrm N_p(0,\Sigma)\in\mathscr N(\D)$ under the mappings
    \begin{align}
 \label{eq:tpd} \left((D,L)\mapsto LD^{-1}L^{\top}\right)&:\Theta_{\D}\rightarrow \mathrm{P}_{\D}\\
 \label{eq:tpdd}\left((D,L)\mapsto \left(LD^{-1}L^{\top}\right)^{-1}\right)&:\Theta_{\D}\rightarrow \mathrm{PD}_{\D}.  
  \end{align}
  
%  We will call these distributions  the DAG-Wishart and the inverse DAG-Wishart distributions on $\mathrm{P}_{\D}$ and $\mathrm{PD}_{\D}$, respectively.
  
   \noindent Viewed inside the full space $\mathrm{S}_p(\mathbb{R})$, every
lower-dimensional DAG model has ambient Lebesgue measure zero.  Its intrinsic
geometry is more informative: for a perfect DAG, the precision model is open
in $\mathrm{Z}_{\D}$, whereas for a non-perfect DAG it is a curved subset of
the larger moral-graph space $\mathrm{Z}_{\D^m}$.  Analogous intrinsic-coordinate
issues arise for covariance matrices.  One can derive a precision density with
respect to Hausdorff measure, but even for small non-perfect DAGs its area
factor complicates posterior analysis; Supplemental section C gives the
precise construction and an example.
 
 \indent To overcome this problem, we follow what was done for the hyper inverse Wishart in \cite{lauritzen96*} or for the type I Wishart in \cite{letac07*} and we work with the projections of $\mathrm{P}_{\D}$ and $\mathrm{PD}_{\D}$ onto the Euclidean space that only retain the functionally independent elements of the precision and covariance matrices of Gaussian DAG models. 
 
  The projected spaces, as we shall see, are subsets of incomplete matrices, which we call  the incomplete precision  space $\Rd$ and the incomplete covariance space $\Sd$, respectively. The precise definitions are as follows.
 \begin{definition}\label{def:incomplete} Let $\D=(V, E)$ be a DAG\footnote{Note an important convention here that the edge set $E$ contains all the loops (see Supplemental section A for details).} and $\D^{u}=(V, E^{u})$ its undirected version.
  \begin{itemize}
  \renewcommand{\labelitemi}{$(a)$}
  \item Let $\mathrm{Z}_{\D}\subseteq \R^{p\times p}$ denote the real
  linear space of $p\times p$ symmetric matrices $A$ such that, for
  $i\ne j$, $A_{ij}=0$ whenever neither $(i,j)$ nor $(j,i)$ belongs to
  $E$.  Thus the allowed off-diagonal positions are the edges of the
  skeleton $\D^u$.  Because $E$ includes the loops, the dimension of
  $\mathrm{Z}_{\D}$ is $|E|$.
  \renewcommand{\labelitemi}{$(b)$}
  \item Let $\mathrm{I}_{\D}$ denote the real linear space of  symmetric functions $\Gamma=\left(\Gamma_{ij}\right)_{(i,j)\in E^{u}}$, i.e., $\Gamma_{ij}=\Gamma_{ji}\in\R$ for each $(i,j)\in E$. An element $\Gamma\in \mathrm{I}_{\D}$  is called a (symmetric) $\D$-incomplete matrix, and can be considered as a matrix in $\mathrm{S}_{p}(\R)$ where only the entries corresponding to the edges of $\D$ are specified and the rest are unspecified. The projection mapping from $\mathrm{S}_{p}(\R)$ onto $\mathrm{I}_{\D}$ is denoted by $A\mapsto A^E$ 
  \renewcommand{\labelitemi}{$(c)$}
  \item For $\Gamma\in\mathrm{I}_{\D}$ let $\left(\Gamma\right)^0$ denote the $|V|\times |V|$ matrix
  \[
  \left(\Gamma\right)^0_{ij}=\begin{cases} \Gamma_{ij}& \text{ if $(i,j)\in E^{u}$},\\ 0& \text{otherwise.}\end{cases}
  \]
  \end{itemize}
    \noindent Note that $\left(\Gamma\right)^0$ fills or completes the unspecified positions with zeros to obtain a full matrix in $\mathrm{Z}_{\D}$. For each clique $c$ of  $\D$  the restriction of  $\Gamma$  on $c$, denoted by $\Gamma_c$,  is a full matrix. Moreover, $\Gamma$ is uniquely determined by the blocks of matrices $(\Gamma_c:   c\in \mathscr{C}_{\D})$, where $\mathscr{C}_{\D}$ denotes the set of cliques of $\D^u$.
    \begin{itemize}
    \renewcommand{\labelitemi}{$(d)$}
  \item Let $\mathrm{Q}_{\D} \subseteq \mathrm{I}_{\D}$ denote the set of $\D$-incomplete matrices $\Gamma\in \mathrm{I}_{\D}$ such that  $\Gamma_{c}$ is positive definite for each clique $c\in \mathscr{C}_{\D}$. Each element of $\mathrm{Q}_{\D}$ is said to be a partially positive definite matrix over $\D$. 
   \renewcommand{\labelitemi}{$(e)$}
    \item Let  $\mathcal{A}\subseteq\mathrm{S}_{p}(\R)$. We say that a  $\D$-incomplete matrix $\Gamma\in\mathrm{I}_{\D}$ can be completed in $\mathcal{A}$ if there exists a matrix  $A\in \mathcal{A}$  such that   $A_{ij}=\Gamma_{ij}$  for each  $(i,j)\in E$, i.e., $A^E=\Gamma$. We refer to $A$ as a completion of $\Gamma$ in $\mathcal{A}$.
		 \renewcommand{\labelitemi}{$(f)$}
		\item The space of incomplete precision matrices over $\D$, denoted by $\Rd$, is the set of $\Gamma\in \mathrm{I}_{\D}$ that can be completed in the space of precision matrices $\mathrm{P}_{\D}$. 
		\renewcommand{\labelitemi}{$(g)$}
		\item The space of incomplete covariance matrices over $\D$, denoted by $\Sd$, is the set of $\Upsilon\in \mathrm{I}_{\D}$ that can be completed in the space of covariance matrices $\mathrm{PD}_{\D}$. 
    \end{itemize}
    \end{definition}
  \begin{Rem}
  If $\mathcal{A}$ is the set of positive definite matrices $\mathrm{PD}_{p}(\R)$, then the completion in $\mathcal{A}$ reduces to the standard definition of positive definite completion \cite{Grone84}. We shall consider below the positive definite completion of partially positive precision/covariance  matrices that correspond to DAGs (vs. those that correspond to undirected graphs as in \cite{Grone84}). Note that an incomplete matrix $\Gamma\in\rm{I}_{\D}$ has a positive definite completion only if $\Gamma\in \mathrm{Q}_{\D}$, i.e., it is partially positive definite over $\D$. 
  \end{Rem}
  %
  % % % % % % % % % % % % % % % % % % % %
  %%%%%%%%%%%%%%%%%%%%%%%%%%%%%%%%%%%%%%%%%
  % % % % % % % % % % % % % % % %
  % % % % % % % % % % % % % % % % % % % % % % %
  % % % % % % % % % % % % % % % % % % % % % % % % %
  % % % % % % % % % % % % % % % % % % % % % % % % % % % % %
  
\subsection{The space of incomplete precision matrices}\label{subs:Rd}
We first recall the definition of $\Rd$ in Definition \ref{def:incomplete} and the following result from \cite{benraj12}.
 
      \begin{proposition}\cite{benraj12}\label{prop:completion_in_PD}
Let $\Gamma$ be a $\D$-partial matrix in $\mathrm{I}_{\D}$. If \ $\Gamma_{11}\neq 0$, then
\begin{itemize}
\renewcommand{\labelitemi}{$(a)$}
\item Almost everywhere (with respect to  the Lebesgue measure on $\mathrm{I}_{\D}$), there exist a unique lower triangular matrix $L\in \mathcal{L}_{\D}$  and a unique diagonal matrix $\Lambda \in \R^{p\times p}$ such that 
$\widehat{\Gamma}=L\Lambda L^{\top}$ is a completion of $\Gamma$.
\renewcommand{\labelitemi}{$(b)$}
\item The matrix $\widehat{\Gamma}$ is the unique positive definite completion of  $\Gamma$ in $\mathrm{P}_{\D}$ if and only if  the diagonal entries of $\Lambda$ are all strictly positive.
\end{itemize}
\end{proposition}

Proposition \ref{prop:completion_in_PD} is of interest to us, because it explicitly shows that without loss of generality every precision matrix  $\Omega\in \Pd$ can be represented by a $\D$-incomplete matrix which only consists of the free parameters of $\Omega$, i.e., $\Omega^E$. The rest, entries corresponding to the missing edges of the DAG, can be discarded, as whenever needed they can be obtained from $\Omega^E$ according to a constructive completion procedure given by the proof of Proposition \ref{prop:completion_in_PD}. We re-formalize this as follows.

 \begin{corollary} \label{cor:pg_rg}
 The projection
 $\Omega\mapsto\Omega^E:\mathrm P_{\D}\rightarrow\mathrm R_{\D}$
 is a diffeomorphism, with inverse
 $\Upsilon\mapsto\widehat\Upsilon$.  In particular,
 $\mathrm R_{\D}$ is open in
 $\mathrm I_{\D}\cong\R^{|E|}$.
 \end{corollary} 

 \begin{proof}
 The projection is linear.  Proposition~\ref{prop:completion_in_PD}
 gives a unique inverse on $\mathrm R_{\D}$.  Recovering the modified
 Cholesky pivots and free entries, and then the omitted precision
 entries, uses rational operations whose pivot denominators are
 strictly positive on $\mathrm R_{\D}$.  Thus the inverse is smooth.
 Equivalently, the coordinate map in \eqref{eq:rd} is bijective and has
 the nonzero Jacobian in Lemma~\ref{lem:jacobian_of_psi}; the inverse
 function theorem gives openness and a smooth local inverse everywhere.
 Uniqueness makes those local inverses the global completion map.
 \end{proof}
 
 % % % % % % % % % % % % % %

  \subsection{The DAG-Wishart distribution on $\Rd$}\label{sub:Rd}
 In light of Corollary \ref{cor:pg_rg} we identify $\mathrm{P}_{\D}$ with $R_{\D}$  through the bijection $\Omega\mapsto \Omega^{E}$. Note that $\Rd$, unlike $\mathrm{P}_{\D}$, is open in its affine support $\mathrm{I}_{\D}$ and, as a consequence of Corollary \ref{cor:pg_rg}, diffeomorphic to $\Theta_{\D}$. Recall that we refer to $\mathrm{R}_{\D}$ as the space of incomplete precision matrices over $\D$. Now let $\pi_{U,\alpha}^{\mathrm{R}_{\D}}$ denote the image  of  $\pi^{{\Theta}_{\D}}_{U,\alpha}$ under the mapping
 % %
\begin{equation}\label{eq:rd}
 \psi\equiv \left(\left(L,D\right)\mapsto\left(LD^{-1}L^{t}\right)^{E}\right):\Theta_{\D}\rightarrow\mathrm{R}_{\D}
\end{equation}
% % %
 \noindent Since  $\mathrm{R}_{\D}$  is  an open subset of the Euclidean space  $\R^{|E|}$, the distribution $\pi_{U,\alpha}^{\mathrm{R}_{\D}}$ has a density with respect to  the Lebesgue measure on $\mathrm{R}_{\D}$. Hence, in light of the diffeomorphism $\Omega\mapsto \Omega^{E}$, in both a natural and practical sense, we define $\pi_{U,\alpha}^{\mathrm{R}_{\D}}$ as the DAG-Wishart distribution on the space of incomplete precision matrices $\mathrm{R}_{\D}$. To derive the density of $\pi_{U,\alpha}^{\mathrm{R}_{\D}}$ we need to compute the Jacobian of the mapping $\psi$ in \eqref{eq:rd}. The  Jacobian  of $\psi$ is a variant of similar transformations found in \cite{roverato00*, khare09a*}. For completeness we still compute this Jacobian  in the following lemma. The proof is given in Supplemental section B subsection 2.7. 
 \begin{lemma}\cite{roverato00*, khare09a*} \label{lem:jacobian_of_psi}
 The Jacobian of the mapping $\psi:(D, L)\mapsto\left( LD^{-1}L^{\top}\right)^{E}$ is $\prod_{j=1}^{p} D_{jj} ^{-(pa_j+2)}$.
 \end{lemma}
 We now proceed to express the density of $\pi_{U,\alpha}^{\mathrm{R}_{\D}}$ and some of its properties. The  proofs are immediate results of Lemma \ref{lem:jacobian_of_psi} and the iterative construction of $\pi^{\Theta_{\cal D}}_{U,\alpha}$. 
 \begin{theorem}\label{thmn:rd}
 Let $\Upsilon$ be the image of $\left(L,D\right)\sim \pi_{U,\alpha}^{\Theta_{\D}}$ under the mapping $\psi$. Then
 \begin{itemize}
 \renewcommand{\labelitemi}{a)}
 \item The density of $\Upsilon\sim \pi^{\mathrm{R}_{\D}}_{U, \alpha}$  with respect to  the standard Lebesgue measure on  $\mathrm{R}_{\D}$ is given by
 \[
 z_{\D}(U,\alpha)^{-1}\exp\left\{-\frac{1}{2}\tr(\widehat{\Upsilon} U)\right\}\prod_{i=1}^{p} D_{ii}^ {-\frac{1}{2}\alpha_i+pa_i+2},
 \]
 where  $D_{ii}=\left(\widehat{\Upsilon}^{-1}\right)_{ii|\prec i\succ}$ is explicitly a function of $\Upsilon$ and $z_{\D}$ is defined in \eqref{eq:zd}.
 \renewcommand{\labelitemi}{b)}
 \item For every symmetric \(K\) such that \(U+2K\succ0\), the
 completion-based Laplace functional is
 \[
 E\{\exp[-\tr(K\Omega)]\}
 =\frac{z_{\D}(U+2K,\alpha)}{z_{\D}(U,\alpha)}.
 \]
 To state the ordinary Euclidean transform without double counting, define
 \[
 \langle H,\Upsilon\rangle_E
 =\sum_i H_{ii}\Upsilon_{ii}
  +\sum_{\substack{i>j\\(i,j)\in E}}H_{ij}\Upsilon_{ij},
 \]
 and let \(K_H\) have diagonal entries \((K_H)_{ii}=H_{ii}\), edge
 entries \((K_H)_{ij}=(K_H)_{ji}=H_{ij}/2\), and zeros elsewhere.  If
 \(U+2K_H\succ0\), then
 \[
 E\{\exp[-\langle H,\Upsilon\rangle_E]\}
 =\frac{z_{\D}(U+2K_H,\alpha)}{z_{\D}(U,\alpha)}.
 \]
 The factors of one half compensate for the two appearances of each
 off-diagonal entry in a symmetric trace.  A general \(K\) in the first
 display can additionally evaluate completed entries.
 \renewcommand{\labelitemi}{c)}
 \item The full precision-matrix mean, not merely the mean of its
 incomplete coordinates, is
 \begin{equation}\label{eq:full-precision-mean}
 \mathbb{E}(\Omega)
 =\sum_{j=1}^{p}(\alpha_j-pa_j-2)
   \left(U_{\preceq j\succeq}^{-1}\right)^0
 -\sum_{j=1}^{p}(\alpha_j-pa_j-3)
   \left(U_{\prec j\succ}^{-1}\right)^0.
 \end{equation}
 Consequently, \(\mathbb{E}(\Upsilon)=\mathbb{E}(\Omega)^E\).
 \end{itemize}
 % % %
 \end{theorem}

\begin{proof}
Part (a) is the change-of-variables formula using
Lemma~\ref{lem:jacobian_of_psi}; part (b) follows by combining the
Laplace kernel with the parameter \(U\).  For part (c), write
\(\Omega=\sum_jD_{jj}^{-1}v_jv_j^\top\), where \(v_j\) has entry one at
\(j\), entries \(\ell_j=-\beta_j\) at \(\mathrm{pa}(j)\), and zero
elsewhere.  Proposition~\ref{prop:nodewise} gives
\[
E(D_{jj}^{-1})=
\frac{\alpha_j-pa_j-2}{U_{jj\mid\prec j\succ}}.
\]
Taking the first two conditional moments of \(\ell_j\) and applying the
block-inverse identity for \(U_{\preceq j\succeq}^{-1}\) yields the
\(j\)th summand in \eqref{eq:full-precision-mean}.  Summing over
vertices proves the claim.
\end{proof}
 %%%%%%%%& % % % % % % % % % % % %
 
 \begin{figure}[H]
   \centering
   \begin{tikzpicture}[x=1cm,y=1cm]
     \coordinate (cycleone) at (-1.05,-0.75);
     \coordinate (cycletwo) at (-1.05,0.75);
     \coordinate (cyclefour) at (1.05,0.75);
     \coordinate (cyclethree) at (1.05,-0.75);
     \foreach \from/\to in {two/one,three/one,four/two,four/three}
       \draw[dag edge,shorten <=2.5mm,shorten >=2.8mm]
         (cycle\from) -- (cycle\to);
     \node[dag vertex] at (cycleone) {$1$};
     \node[dag vertex] at (cycletwo) {$2$};
     \node[dag vertex] at (cyclethree) {$3$};
     \node[dag vertex] at (cyclefour) {$4$};
   \end{tikzpicture}
   \caption{A four-node DAG whose skeleton is a 4-cycle.}
   \label{fig-4}
 \end{figure}
%%%%

 % % % % % % % % % % %
 \begin{Ex}
 Let $\D$ be the DAG given by Figure \ref{fig-4}. Then the DAG-Wishart density $\pi_{U,\alpha}^{\mathrm{R}_{\D}}$ is given by
 \[
 \scalebox{.85}{$\displaystyle
 \pi_{U,\alpha}^{\mathrm{R}_{\D}}\left(\Upsilon\right)=
 z_{\D}(U,\alpha)^{-1}\exp\left\{-\frac{1}{2}\tr\left(\widehat{\Upsilon}U \right)\right\}D_{11}\left(\Upsilon\right)^ {-\frac{1}{2}\alpha_1+4}D_{22}\left(\Upsilon\right)^ {-\frac{1}{2}\alpha_2+3}D_{33}\left(\Upsilon\right)^ {-\frac{1}{2}\alpha_3+3}D_{44}\left(\Upsilon\right)^ {-\frac{1}{2}\alpha_4+2}$,}
 \]
 where, using Proposition \ref{prop:completion_in_PD}, $\widehat{\Upsilon}$ and $D_{ii}\left(\Upsilon\right)$ are computed as follows. 
 \[
 \scalebox{.83}{$\displaystyle
\widehat{\Upsilon}=
 \begin{pmatrix}
 \Upsilon_{11} & \Upsilon_{21} & \Upsilon_{31} & 0 \\
 \Upsilon_{21} & \Upsilon_{22} & \dfrac{\Upsilon_{21}\Upsilon_{31}}{\Upsilon_{11}} & \Upsilon_{42} \\
 \Upsilon_{31} & \dfrac{\Upsilon_{21}\Upsilon_{31}}{\Upsilon_{11}}& \Upsilon_{33} & \Upsilon_{43} \\
 0 & \Upsilon_{42} & \Upsilon_{43} & \Upsilon_{44}
 \end{pmatrix}
 $}
 \]
 \[
\begin{aligned}
 D_{11}&=\Upsilon_{11}^{-1},\\
 D_{22}&=\left(\Upsilon_{22}-\frac{\Upsilon_{21}^{2}}{\Upsilon_{11}}\right)^{-1},
 &D_{33}&=\left(\Upsilon_{33}-\frac{\Upsilon_{31}^{2}}{\Upsilon_{11}}\right)^{-1},\\
 D_{44}&=\left(
 \widehat\Upsilon_{44}
 -\widehat\Upsilon_{4,\{1,2,3\}}
  \widehat\Upsilon_{\{1,2,3\}}^{-1}
  \widehat\Upsilon_{\{1,2,3\},4}
 \right)^{-1}.
\end{aligned}
\]
 \end{Ex}
 \medskip

% % % % % % % % % % % % % % % % % % % % % %
\section{The inverse DAG-Wishart distribution on the space of incomplete covariance matrices}\label{sec:Sd}
We now derive the covariance-coordinate image of the DAG--Wishart,
which parallels the hyper-inverse Wishart and the inverse type-II
Wishart \(\mathrm{IW}_{\Pg}\).  We first define incomplete covariance
matrices and recall the completion results from \cite{benraj12} needed
to construct this image.

%%%%%%%%%%%%%%%%%%%%%%%%%%%%%%%%%%%%%%%%%%%%%%%%%%%%%%%
% % % % % % % % % % % % % % %
% % % % % % % % % % % % % % % % % %
% % % % % % % % % % % % % % % % % % % % %
% % % % % % % % % % % % % % % % % % % % % % % % % %
% % % % % % % % % % % % % % % % % % % % % % % % % % % %

\subsection{The space of incomplete covariance matrices}\label{sub:completion_in_pdd} 

Recall that $\mathrm{PD}_{\D}$  is the space of covariance matrices for the Gaussian DAG model $\mathscr{N}(\D)$, the elements of which, according to \eqref{eq:ndmp}, can be characterized as:
\begin{equation}\label{eq:sigma}
\Sigma\in \mathrm{PD}_{\D}\Longleftrightarrow  \Sigma\succ 0\: \: \& \:\: \Sigma_{[i\nsucc}=\Sigma_{[i\succ}\Sigma_{\prec i\succ}^{-1}\Sigma_{\prec i\nsucc},\: \text{for each $i\in V$}. 
\end{equation}
The above characterization allows us to identify $\mathrm{PD}_{\D}$ with the functionally independent elements of $\Sigma$. The following proposition is a key ingredient in this identification. 
\begin{proposition}\cite{benraj12}\label{prop:completion_in_PDD}
Let $\Gamma\in\mathrm{Q}_{\D}$, then
\begin{enumerate}
\item There exists a completion process of polynomial complexity that can determine whether $\Gamma$ can be completed in $\mathrm{PD}_{\D}$;
\item If a completion exists, this completion is unique and can be determined constructively using the following process:
\begin{itemize}
\renewcommand{\labelitemi}{$i)$}
\item Set $\Sigma_{ij}=\Gamma_{ij}$ \; for each $(i,j)\in E$  and set  $j=p$.
\renewcommand{\labelitemi}{$ii)$}
\item If  $j>1$, then set $j=j-1$ and proceed to the next step, otherwise $\Sigma$ is successfully completed.
\renewcommand{\labelitemi}{$iii)$}
\item If $\Sigma_{\preceq j\succeq}> 0$, then proceed\footnote{Note that for each $j$, the submatrix $\Sigma_{\preceq j\succeq}$ is fully determined by step {\it(ii)}} to the next step, otherwise the completion in $\mathrm{PD}_{\D}$  does not exist.
\renewcommand{\labelitemi}{$iv)$}
\item If $\Sigma_{\nprec j]}$ is non-empty, then set $\Sigma_{\nprec j]}=\Sigma_{\nprec j\succ}\Sigma_{\prec j\succ}^{-1}\Sigma_{\prec j]},$ $\Sigma_{[ j\nsucc}=\Sigma_{\nprec j ]}^{t}$ and return to step $(2).$
\end{itemize}
\end{enumerate}
\end{proposition}

\begin{Rem}
Note once more that the procedure in Proposition \ref{prop:completion_in_PDD} itself determines if $\Gamma$ can be completed in $\mathrm{PD}_{\D}$. It is clear from Step {\it(iii)} above that the necessary and sufficient condition for the existence of a positive definite completion is that, for each $j\in V$, the covariance sub-matrix $\Sigma_{\preceq j\succeq}> 0$ and not just  $\Sigma_{\prec j \succ}>0$. Furthermore, the completion procedure in Proposition \ref{prop:completion_in_PDD} can terminate midway.
\end{Rem}
From Definition \ref{def:incomplete} recall that $\mathrm{S}_{\D}$ denotes the set of $\Gamma\in \mathrm{I}_{\D}$ that can be completed in $\mathrm{PD}_{\D}$. We call this set the space of incomplete covariance matrices over $\D$. The next corollary formalizes the fact that $\mathrm{S}_{\D}$ can be identified with $\mathrm{PD}_{\D}$.
Its proof is immediate from Proposition  \ref{prop:completion_in_PDD} above. 
\begin{corollary}\label{cor:pdd_sd}
The projection
$\Sigma\mapsto\Sigma^E:\mathrm{PD}_{\D}\rightarrow\mathrm S_{\D}$
is a diffeomorphism with inverse
$\Gamma\mapsto\widetilde\Gamma$, where $\widetilde\Gamma$ is the
completion constructed in Proposition~\ref{prop:completion_in_PDD}.
In particular, $\mathrm S_{\D}$ is open in
$\mathrm I_{\D}\cong\R^{|E|}$.
\end{corollary}
\begin{proof}
The projection is linear and Proposition~\ref{prop:completion_in_PDD}
gives a unique inverse.  At step $j$, the completion uses
$\Sigma_{\prec j\succ}^{-1}$; its determinant is positive throughout
$\mathrm S_{\D}$.  Every completed entry is therefore a rational,
hence smooth, function of the incomplete coordinates.  The Jacobian in
Lemma~\ref{lem:inverse_mapping} is strictly positive, so the inverse
function theorem also gives openness and smoothness of the inverse.
\end{proof}
\begin{Rem}\label{rem:QG_SG}
Suppose $\D$ is perfect. Then $\mathrm{PD}_{\D}$ is identical to $\mathrm{PD}_{\D^{\mathrm{u}}}$ and, therefore, by the completion result in Grone et al. \cite{Grone84}, every incomplete matrix in $\mathrm{Q}_{\D}$  can be completed in $\mathrm{PD}_{\D}$. Hence for $\D$ perfect, $\mathrm{S}_{\D}$ and $\mathrm{Q}_{\D}$ are identical.
\end{Rem}
%

%%%%%%%%%%%%%%%%%%%%%%%%%%%%%%%%%%%%%%%%%%%%%%%%%%%%%%%%
% % % % % % % % % % % % % % % % % % % % % % % % % % % %
% % % % % % % % % % % % % % % % % % % % % % % % % % % % % % % %
% % % % % % % % % % % % % % % % % % % % % % % % % % % % % % % % %
% % % % % % % % % % % % % % % % % % % % % % % % % % % % % % % % % % %
% % % % % % % % % % % % % % % % % % % % % % % % % % % % % % % % % % % % % % %

\subsection{The inverse DAG-Wishart distribution on $\mathrm{S}_{\D}$}\label{sub:wsg}

Let $\pi_{U,\alpha}^{\mathrm{S}_{\D}}$ denote the image of
$\pi_{U,\alpha}^{\Theta_{\D}}$ under
\[
(D,L)\longmapsto (L^{-\top}DL^{-1})^{E}:
\Theta_{\D}\longrightarrow \mathrm{S}_{\D}.
\]
In parallel to our notation $\pi_{U,\alpha}^{\mathrm{R}_{\D}}$, we call
this the inverse DAG--Wishart distribution on incomplete covariance
coordinates.  To derive its Lebesgue density, we first compute the
Jacobian of
$(\Sigma^{E}\mapsto\Sigma^{-E}):\mathrm{S}_{\D}\rightarrow
\mathrm{R}_{\D}$, where
\(\Sigma^{-E}=(\Sigma^{-1})^E\) and \(\Sigma\) is the unique completion
of \(\Sigma^E\) in \(\mathrm{PD}_{\D}\).

\begin{lemma}\label{lem:inverse_mapping}
Let $\D=(V,E)$  be an arbitrary DAG, then the Jacobian of the mapping $ (\Sigma^{-E}\mapsto \Sigma^{E}): \mathrm{R}_{\D}\rightarrow \mathrm{S}_{\D}$ is given by $\prod_{i=1}^{p}\dfrac{\det\Sigma_{\preceq i\succeq}^{pa_i+2}}{\det\Sigma_{\prec i\succ}^{pa_i+1}}$.
\end{lemma}

\begin{proof}
Write the mapping $\Sigma^{-E}\mapsto\Sigma^E$ as the composition
\begin{align*}
&(\Sigma^{-E}\mapsto\times_{i=1}^{p} (\Sigma_{ii|\prec i\succ}, \Sigma_{\prec i\succ}^{-1}\Sigma_{\prec i]}):\mathrm{R}_{\D}\rightarrow \Xi_{\D};\\
&(\times_{i=1}^{p}(\Sigma_{ii|\prec i\succ}, \Sigma_{\prec i\succ}^{-1}\Sigma_{\prec i]})\mapsto \Sigma^{E}):\Xi_{\D}\rightarrow \mathrm{S}_{\D}.
\end{align*}
The first map is the inverse of the free-coordinate version of $\psi$
in Lemma~\ref{lem:jacobian_of_psi}; its Jacobian is
\[
\prod_{i=1}^p\lambda_i^{pa_i+2},
\qquad \lambda_i=\Sigma_{ii\mid\prec i\succ}.
\]

For the second map, it is simplest to calculate the inverse
$\Sigma^E\mapsto(\lambda_i,\beta_i)_{i=1}^p$.  Remove vertex~1 and
apply induction to the ancestral subgraph on
$\{2,\ldots,p\}$.  Conditional on that subgraph, the remaining block
map is
\[
(\Sigma_{11},\Sigma_{\prec1]})
\longmapsto
\left(\Sigma_{11\mid\prec1\succ},
      \Sigma_{\prec1\succ}^{-1}\Sigma_{\prec1]}\right),
\]
whose block-triangular derivative has determinant
$\det(\Sigma_{\prec1\succ})^{-1}$.  Hence
\[
J\{\Sigma^E\mapsto(\lambda_i,\beta_i)_{i=1}^p\}
=\prod_{i=1}^p\det(\Sigma_{\prec i\succ})^{-1},
\]
and the forward map from $\Xi_{\D}$ to $\mathrm S_{\D}$ has Jacobian
$\prod_i\det(\Sigma_{\prec i\succ})$.  Multiplying the two forward
Jacobians and using
\(
\lambda_i=\det(\Sigma_{\preceq i\succeq})/
             \det(\Sigma_{\prec i\succ})
\)
gives
\[
\prod_{i=1}^{p}
\frac{\det(\Sigma_{\preceq i\succeq})^{pa_i+2}}
     {\det(\Sigma_{\prec i\succ})^{pa_i+1}},
\]
as claimed.
\end{proof}

We now proceed to state the functional form of the density of $\pi_{U,\alpha}^{\mathrm{S}_{\D}}$ with respect to Lebesgue measure.
\begin{corollary}\label{cor:density_SG}
Let  $\Sigma\sim \pi_{U,\alpha}^{\mathrm{PD}_{\D}}$ and let $\Gamma=\Sigma^{E}$, i.e., $\Sigma$ is the completion of $\Gamma$ in $\mathrm{PD}_{\D}$. Then the  density of  $\Gamma\sim \pi_{U,\alpha}^{\mathrm{S}_{\D}}$ with respect to  Lebesgue measure is given by
\begin{equation}\label{eq:density_SG}
 z_{\D}(U,\alpha)^{-1}\exp\left\{-\frac{1}{2}\tr(\Sigma^{-1} U)\right\}\prod_{i=1}^{p}\dfrac{\det \Sigma_{\preceq i\succeq}^ {-\frac{1}{2}\alpha_i}}{\det\Sigma_{\prec i\succ}^{-\frac{1}{2}\alpha_i+1}}.
\end{equation}
\end{corollary}
\begin{figure}[htbp]
  \centering
  \ColliderDAG
  \caption{The collider DAG studied in Example \ref{ex:Sd}.}
  \label{fig-33b}
\end{figure}

\begin{Ex}\label{ex:Sd}
Consider the  DAG $\D$ given in Figure \ref{fig-33b}. Then the inverse DAG-Wishart on $\D$ is given by
\[
\pi_{U,\alpha}^{\mathrm{S}_{\D}}(\Gamma)=z_{\D}(U,\alpha)^{-1}\exp\left\{-\frac{1}{2}\tr(\Sigma^{-1} U)\right\} D^ {-\frac{1}{2}\alpha_1}_{11}D_{22}^ {-\frac{1}{2}\alpha_2}D_{33}^ {-\frac{1}{2}\alpha_3}\det(\Sigma_{\prec 1\succ})^{-1},
\]
where $\Sigma$, the completion of $\Gamma$, is simply computed as
\[
\Sigma=
\begin{pmatrix}
\Gamma_{11}&\Gamma_{12}&\Gamma_{13}\\
\Gamma_{21}&\Gamma_{22}&0\\
\Gamma_{31}&0&\Gamma_{33}
\end{pmatrix}.
\]
\end{Ex}
\begin{Rem}
For a decomposable graph, compatible perfect-DAG specializations of
$\pi^{\Sd}_{U,\alpha}$ have the same block conditional structure as
the inverse type-II Wishart of \cite{letac07*}.  This statement does not
identify the unrestricted vertexwise family with the entire type-II
family; their shape parameterizations and admissibility constraints
differ.  For the homogeneous covariance-graph subclass, the
specialization described in Supplemental Section~B agrees with the
inverse Wishart construction of Khare and Rajaratnam
\cite{khare09a*}.
\end{Rem}
%%%%%%%%%%%%%%%%%%%%%%%%%%%%
% % % % % % % % % % % % % % %
% % % % % % % % % % % % % % % % %
% % % % % % % % % % % % % % % % % % % % %
% % % % % % % % % % % % % % % % % % % % % % % % % %
% % % % % % % % % % % % % % % % % % % % % % % % % % % % %
% % % % % % % % % % % % % % % % % % % % % % % % % % % % % %
\subsection{Properties of the inverse DAG-Wishart distributions}\label{sub:prop_wsg}
One of the main useful properties of the inverse DAG--Wishart
$\pi_{U,\alpha}^{\mathrm S_{\D}}$ is its nodewise directed
hyper-Markov factorization.  For an arbitrary DAG this follows directly
from Proposition~\ref{prop:nodewise} under the diffeomorphism in
Corollary~\ref{cor:pdd_sd}; it is not inferred from the decomposable-graph
Theorem~4.4 of \cite{letac07*}.  For a perfect DAG, the two results have
parallel block conditional interpretations.
\begin{theorem}\label{th:strong} If ~$\Sigma^{E}\sim {\pi}_{U,\alpha}^{\mathrm{S}_{\D}}$, then \\
 $i)$~ $\left\{(\Sigma_{ii|\prec i\succ}, \Sigma_{\prec i\succ}^{-1}\Sigma_{\prec i]}: i\in V\right\}$ are mutually independent; thus ${\pi}^{\mathrm{S}_{\D}}_{U,\alpha}$ has the stated strong directed hyper-Markov factorization.\\
 $ii)$~ The distribution of  $\Sigma_{ii|\prec i\succ}$  and  $\Sigma_{\prec i\succ}^{-1}\Sigma_{\prec i]}|\Sigma_{ii|\prec i\succ}$ are, respectively, given by
 \begin{equation}\label{hypersg1}
\Sigma_{ii|\prec i\succ}\sim IG(\frac{\alpha_i}{2}-\frac{pa_i}{2}-1,  \frac{1}{2}U_{ii|\prec i\succ}) ,~\text{and}
\end{equation}
\begin{equation}\label{hypersg2}
\Sigma_{\prec i\succ}^{-1}\Sigma_{\prec i]}|\Sigma_{ii|\prec i\succ}\sim \rm{N}_{pa_i}(U^{-1}_{\prec i\succ}U_{\prec i]}, \Sigma_{ii|\prec i\succ}U^{-1}_{\prec i\succ} ).
\end{equation}
 \end{theorem}

\begin{proof}
Under the covariance map,
$D_{ii}=\Sigma_{ii\mid\prec i\succ}$ and
$\beta_i=\Sigma_{\prec i\succ}^{-1}\Sigma_{\prec i]}=-L_{\prec i]}$.
The mutual independence and the two displayed laws are therefore
exactly Proposition~\ref{prop:nodewise} expressed in covariance
coordinates.
\end{proof}

The next result computes the expectation of the \emph{full} covariance
matrix, including entries that are not retained in \(\Sigma^E\).

\begin{proposition}\label{exofsigma}
Suppose \(\Sigma=L^{-\top}DL^{-1}\) is induced by
\(\pi^{\Theta_{\D}}_{U,\alpha}\), and assume
\(\alpha_i>pa_i+4\) for every \(i\).  Define
\[
m_i=U_{\prec i\succ}^{-1}U_{\prec i]},
\qquad
r_i=\frac{U_{ii\mid\prec i\succ}}{\alpha_i-pa_i-4},
\qquad
B_i=r_iU_{\prec i\succ}^{-1}+m_im_i^\top .
\]
For an empty parent set, take \(m_i\) and \(B_i\) to be empty and
\(r_i=U_{ii}/(\alpha_i-4)\).  The matrix
\(M=\mathbb{E}(\Sigma)\) is obtained for \(i=p,p-1,\ldots,1\) as
follows:
\begin{align}
M_{ji}=M_{ij}
  &=M_{j,\prec i\succ}m_i,
  &&j=i+1,\ldots,p,\label{eq:cov-mean-offdiag}\\
M_{ii}
  &=r_i+\tr\!\left(M_{\prec i\succ}B_i\right).
  \label{eq:cov-mean-diag}
\end{align}
When \(\mathrm{pa}(i)=\emptyset\), the first expression is zero and the
trace term is absent.  In particular,
\(\mathbb{E}(\Sigma^E)=M^E\).
\end{proposition}

\begin{proof}
The recursive Gaussian representation associated with the modified
Cholesky factor is
\[
X_i=\beta_i^\top X_{\mathrm{pa}(i)}+\varepsilon_i,
\qquad
\varepsilon_i\sim\mathrm{N}(0,D_{ii}),
\]
with \(\varepsilon_i\) independent of the higher-index variables.
The parameters \((D_{ii},\beta_i)\) are independent of the parameters
that determine \(\Sigma_{\{i+1,\ldots,p\}}\).  Proposition
\ref{prop:nodewise} gives
\[
E(D_{ii})=r_i,\qquad
E(\beta_i)=m_i,\qquad
E(\beta_i\beta_i^\top)=B_i.
\]
Taking expectations in
\(\Sigma_{ji}=\Sigma_{j,\mathrm{pa}(i)}\beta_i\) and
\(\Sigma_{ii}=D_{ii}+
\beta_i^\top\Sigma_{\mathrm{pa}(i)}\beta_i\) proves
\eqref{eq:cov-mean-offdiag}--\eqref{eq:cov-mean-diag}.
\end{proof}

\begin{Rem}\label{rem:expectation-completion}
The completion map is nonlinear.  Consequently, the DAG completion of
\(\mathbb{E}(\Sigma^E)\) is generally not \(\mathbb{E}(\Sigma)\);
likewise, completing \(\mathbb{E}(\Omega^E)\) generally does not give
\(\mathbb{E}(\Omega)\).  Equations~\eqref{eq:full-precision-mean} and
\eqref{eq:cov-mean-offdiag}--\eqref{eq:cov-mean-diag} must be used when
a full-matrix posterior mean is required.
\end{Rem}

%%%%%%%%%%%%%%%%%%%%%%%%%%%%%%%%%%%%%%%%%%%%%%%%%%%%%%%%%%%%%%
\section{Simulation study and applications to real data}\label{sec:empirical}
\label{simul}
This section records the numerical studies associated with the original
preprint and states exactly what they can support.  The structure search
is over DAGs compatible with a \emph{fixed, supplied parent ordering};
it is not unrestricted DAG learning and it does not estimate the
ordering.  The public companion archive contains the principal R
functions and some driver scripts, but not the saved simulation outputs,
the Sachs input used for the figure, or the call-center data.  The
tables below are therefore historical point estimates unless explicitly
identified as reproduced.  They illustrate feasibility and motivate
further evaluation; they do not establish a general scaling rate or
uniform risk dominance.

\subsection{Bayesian model selection via DAG-Wishart prior}
Let \(\mathfrak{D}_{\prec}\) be the set of DAGs compatible with the fixed
ordering.  For \(\mathcal{D}\in\mathfrak{D}_{\prec}\),
\[
p(\mathcal{D}\mid X)\propto
p(X\mid\mathcal{D})p(\mathcal{D}).
\]
For observations modeled with known zero mean and
\(T=X^\top X\), conjugacy gives the exact marginal likelihood
\begin{equation}\label{eq:marginal-score}
p(X\mid\mathcal{D})
=(2\pi)^{-np/2}
\frac{z_{\mathcal{D}}(U+T,\alpha+n\mathbf{1})}
     {z_{\mathcal{D}}(U,\alpha)}.
\end{equation}
The original experiments used a uniform prior over
\(\mathfrak{D}_{\prec}\).  Because each of the \(p(p-1)/2\) admissible
edges can then be present or absent independently, this prior has
edge-inclusion probability \(1/2\) and is not sparsity-favoring.  A
transparent alternative is to write
\(e(\mathcal D)=|E(\mathcal D)|-p\) for the number of non-loop edges
and use
\[
p(\mathcal{D}\mid\rho)\propto
\rho^{e(\mathcal{D})}(1-\rho)^{p(p-1)/2-e(\mathcal{D})},
\]
with fixed \(\rho<1/2\) or a stated beta hyperprior.  The historical
numbers below retain the uniform graph prior so that their target is not
changed retrospectively.

The search heuristic combines stochastic shotgun search (SSS)
\cite{SSS} with starting graphs from the LassoDAG regularization path
\cite{LassoDAG}.  It approximates the maximizer of the posterior score;
it does not sample from the exact graph posterior.  The DAG-W procedure
is specified below.
%%%%%%%%%%%%%%%%%%%%%
\begin{algo}[DAG-W]\label{algo:DAG-W}
Assume the following are given: the standardized data matrix $X$, the hyper-parameters $\alpha$, $U$ and the maximum iteration number $M$. Estimate  $N$ models corresponding to different points on the LassoDAG regularization path, labeled as $\mathcal{D}^{(k)}, k=1, \cdots, N$. Then for each $k = 1, 2, \cdots, N$, do the following.
\begin{enumerate}
\item Let $\mathcal{D}_0 = \mathcal{D}^{(k)}$. Until the maximum iteration number $M$ is achieved:
\begin{enumerate}
\item Select $N_1$ graphs that differ from $\mathcal{D}_0$ by one
admissible edge toggle.  Evaluate their log posterior scores
$s_1,\ldots,s_{N_1}$ using \eqref{eq:marginal-score} and the stated
graph prior, and update the retained list $\mathcal{L}^{(k)}$.
\item Sample the next graph from the current candidates with probability
\[
p_i=\frac{\exp\{\gamma(s_i-s_{\max})\}}
{\sum_j\exp\{\gamma(s_j-s_{\max})\}},
\]
where $\gamma$ is an annealing parameter.  Set the sampled graph as the
new $\mathcal{D}_0$.
\item Return to Step 1-(a).
\end{enumerate}
\item Combine the retained lists \(\mathcal{L}^{(k)}\), \(k=1,\ldots,N\).
\item Return the graph with the largest score as the selected model.
\end{enumerate}
\end{algo}

The initial models come from different points on the LassoDAG
regularization path.  In \cite{LassoDAG}, the penalty parameter
\(\tau_i\) for the Lasso problem at node \(i\) is
\begin{equation}
\label{eq:eq1}
\tau_i = 2\frac{Z^*_{\frac{\kappa}{2p(i-1)}}}{\sqrt{n}},
\end{equation}
where \(Z^*_q\) denotes the \((1-q)\)th quantile of the standard normal
distribution; \cite{LassoDAG} recommends \(\kappa=0.1\).  The
historical comparison uses that setting for LassoDAG.

The historical model-selection study used \(U=I\) and
\(\alpha_i=c\,pa_i+b\), with \(c=1\) and \(b=3\).  This satisfies
properness but implies
\[
D_{ii}\sim IG(1/2,1/2),\qquad
\beta_i\ \hbox{has a multivariate \(t\) law with one degree of freedom}.
\]
Thus the local regression prior is Cauchy-tailed and neither
\(E(D_{ii})\) nor \(E(\beta_i)\) exists marginally; the setting should
not be described through prior means.  The study used \(N=16\) initial
states: fifteen values \(\kappa=(k/15)^4p\), \(k=1,\ldots,15\), and the
LassoDAG recommendation \(\kappa=0.1\).  Values \(\kappa>1\) fall
outside the false-positive-control interpretation in \cite{LassoDAG}
and were used only to diversify starting graphs.  The reported search
settings were \(M=100\), \(N_1=30\), and \(\gamma=0.5\), with reduced
search budgets for the largest dimensions.

The driver code uses \texttt{pcalg::randomDAG} with independent
admissible-edge probability \(0.01\) and regression weights sampled
between \(0.2\) and \(0.8\).  It fixes \(n=100\) and considers
\(p=50,100,200,500,1000,1500,2000\).  For \(p\geq500\), the manuscript
reports nine starts and at most 50 search steps per start.  The public
archive does not contain a \(p=2000\) driver or any saved output.  Its
\(p=1000\) and \(p=1500\) drivers call an absent
\texttt{PerfEval2.R}, request 16 starts and 100 search iterations rather
than the reduced manuscript budgets, and use only two replications per
parallel invocation.  Package versions, hardware, elapsed times, Monte
Carlo standard errors, and a complete mapping from scripts to table rows
were not preserved.  These mismatches and omissions prevent exact reproduction of
Table~\ref{tab:sparsity001table} from the archive alone.

There is also a likelihood mismatch in the archived search path.  The
selection routine centers and scales the observations using sample
statistics, after which the marginal-score routine increments \(\alpha\)
by \(n\).  Under the flat unknown-mean analysis in
Remark~\ref{rem:unknown-mean}, centering instead yields the update
\(\alpha+(n-1)\mathbf1\).  Thus the historical graph scores should be
understood as scores from a data-dependent standardized, zero-mean
working likelihood, not as the corrected unknown-mean marginal
likelihood.

Sensitivity and specificity are reported in the historical table.
With approximately one true edge per 100 admissible edges, specificity
near one can still correspond to many false discoveries, so precision
(positive predictive value), MCC, structural Hamming distance, and
uncertainty should also be reported in a new experiment.  For example,
using the \(p=2000\) averages as plug-in rates gives roughly 3,886 true
and 11,082 false positives for DAG-W (precision about \(26\%\)), versus
1,977 true and 6,531 false positives for LassoDAG (precision about
\(23\%\)).  This calculation is only an interpretation of averaged
rates, not a reconstruction of the missing replicate-level results.

\begin{table}[htbp]
\begin{center}
\begin{tabular}{ |c|cc|cc| }
\hline
&\multicolumn{2}{ |c| }{LassoDAG}& \multicolumn{2}{ |c| }{DAG-W} \\
\hline
p& Sensitivity & Specificity & Sensitivity & Specificity \\ \hline
50 & 0.6156 & 1.0000 &  0.7828 & 0.9980\\ \hline
100 & 0.4826 &$ \sim 1$ & 0.7524 & 0.9977\\ \hline
200 & 0.3969 & $\sim 1$ & 0.7405 & 0.9975\\ \hline
500 & 0.2497 & $\sim 1$ & 0.6517 & 0.9982\\ \hline
1000& 0.1748 & 0.9991 & 0.4248 & 0.9971\\ \hline
1500 & 0.1226 & 0.9981 & 0.2672 & 0.9962\\ \hline
2000 & 0.0989 & 0.9967 & 0.1944 & 0.9944\\ \hline
\end{tabular}
\end{center}
\caption{Historically reported average sensitivity and specificity for
\(n=100\) and admissible-edge probability \(0.01\).  Replicate-level
outputs and Monte Carlo uncertainty are unavailable in the public
archive; the table is not independently reproduced in this revision.}
\label{tab:sparsity001table}
\end{table}

\subsection{Covariance Estimation Performance}
We now consider covariance and precision estimation \emph{conditional
on the true DAG and its ordering}; graph-selection uncertainty is not
included in this experiment.  Following \cite{rajaratnam08*}, the
historical code uses a coordinatewise squared-error loss
\[
L_2(M,\widehat{M}) =
\sum_{(i,j)\in E}(M_{ij}-\widehat{M}_{ij})^2,
\]
where \(E\) includes the diagonal under the convention of
Definition~\ref{def:incomplete}.  This is not the full Frobenius loss.
The full-matrix Stein loss used in the code is
\[ 
L_1(\widehat{M},M) = \tr(\widehat{M}M^{-1}) - \log(\det(\widehat{M}M^{-1}))-p.
\]
Let
\[
\bar\Sigma=E(\Sigma\mid X),\qquad
\bar\Omega=E(\Omega\mid X),
\]
computed from Proposition~\ref{exofsigma} and
\eqref{eq:full-precision-mean} with the posterior parameters.  Their
edge projections are Bayes actions under the coordinatewise squared
loss.  Under the displayed full Stein loss, the Bayes actions are
\[
\delta_\Sigma^{\mathrm{Stein}}=\bar\Omega^{-1},
\qquad
\delta_\Omega^{\mathrm{Stein}}=\bar\Sigma^{-1}.
\]
Neither action is generally DAG-constrained, because expectations and
inversion do not preserve the nonlinear DAG parameter spaces.

The historical precision routine first computes \(\bar\Omega^E\) and
then applies the DAG completion map; denote this graph-constrained
representative by
\(\widehat\Omega_{\mathrm{coord}}
=\mathcal{C}_{\D}(\bar\Omega^E)\).
It has the correct posterior-mean edge coordinates and is therefore
valid for the coordinate loss, but it is not \(\bar\Omega\) and is not
the full-Stein Bayes action.  The third reported estimator transforms
the joint posterior mode in \((D,L)\) coordinates into covariance and
precision matrices.  We denote these plug-ins by
\(\widetilde\Sigma_{\mathrm{chol}}\) and
\(\widetilde\Omega_{\mathrm{chol}}\); they are not called MAP
estimators because modes change under nonlinear reparameterization.
Supplemental Section~D gives the corrected algorithms.

The historical study uses \(\alpha_i=c\,pa_i+3\) and
\(U=uI\), with the values shown in the table, and simulates
\(p=500\) with admissible-edge probability \(0.01\).  The driver
specifies 15 parallel batches of 20 replications per setting, but the
saved outputs are absent, so execution and aggregation cannot be
verified.  Table~\ref{tab:Estimation-p500-Precision-detail} is retained
as a historical report of relative changes from the graph-constrained
MLE.  It contains no Monte Carlo standard errors or intervals, and its
Stein-loss rows evaluate the listed procedures rather than establishing
that each is a Bayes action.

\begin{table}[htbp]
\begin{center}
\begin{tabular}{ |c|c|cc|cc|cc| }
\hline
&&\multicolumn{2}{ |c| }{n=30}& \multicolumn{2}{ |c| }{n=50} &  \multicolumn{2}{ |c| }{n=100} \\
  \hline
  $(c,U)$&Estimator & $L_1$ & $L_2$ & $L_1$ & $L_2$ & $L_1$ & $L_2$\\ 
\hline

\multirow{3}{*}{$(2.5,I(3))$} &$\widehat{\Omega}_{\mathrm{coord}}$            & 41.8\% & 77.9\% & 26.8\% & 56.5\%  & 14.2\% & 29.8\% \\
                                              & $\bar\Sigma^{-1}$     & 45.8\% & 60.2\%  & 29.8\% & 30.7\%  & 15.9\% & 3.6\%\\
                                              & $\widetilde{\Omega}_{\mathrm{chol}}$             & 38.7\% & {\bf 82.0}\%  & 23.9\% & {\bf 63.0\%}  & 12.3\%   &37.9\%\\

\hline

\multirow{3}{*}{$(3,I(3))$} &$\widehat{\Omega}_{\mathrm{coord}}$          & 39.2\% & 80.5\%  & 24.7\% & 60.5\% & 12.9\% & 34.6\% \\
                                              & $\bar\Sigma^{-1}$& 47.4\% & 65.9\%  & 31.1\% & 39.9\% & 16.7\% & 13.8\%\\
                                              & $\widetilde{\Omega}_{\mathrm{chol}}$        & 34.4\% & 81.5\% & 20.1\% & 62.3\% & 9.7\% & 37.7\%\\
\hline
\multirow{3}{*}{$(3.5,I(3))$} &$\widehat{\Omega}_{\mathrm{coord}}$            & 35.9\% & 81.9\% & 21.9\% & 62.8\%  & 11.1\% & 37.4\% \\
                                              & $\bar\Sigma^{-1}$     & {\bf 47.9\%} & 70.1\%  & {\bf 31.6\%} & 47.6\%  & {\bf 17.1\%} & 22.3\%\\
                                              & $\widetilde{\Omega}_{\mathrm{chol}}$             & 29.5\% & 79.9\%  & 15.7\% & 59.7\%  &    6.7\%   & 35.5\%\\
  \hline
  \multirow{3}{*}{$(3,I(2.5))$} &$\widehat{\Omega}_{\mathrm{coord}}$          & 34.5\% & 81.9\% & 20.1\% & {\bf 63.0\%} & 10.3\% & {\bf38.6\%} \\
                                                   & $\bar\Sigma^{-1}$& 47.8\% & 72.4\% & 31.3\% & 51.0\%   &16.8\% & 27.1\%\\
                                                   & $\widetilde{\Omega}_{\mathrm{chol}}$        & 26.6\% & 77.2\% & 13.2\% & 55.9\% & 5.1\% & 32.7\%\\
  \hline
  \multirow{3}{*}{$(3,I(3.5))$} &$\widehat{\Omega}_{\mathrm{coord}}$          & 42.9\% & 77.0\% & 27.0\% & 54.2\% & 14.4\% & 25.8\% \\
                                                   & $\bar\Sigma^{-1}$& 45.6\% & 59.0\% & 29.6\% & 27.3\% &15.7\% & -2.6\%\\
                                                   & $\widetilde{\Omega}_{\mathrm{chol}}$        & 39.6\% & 81.9\% & 24.9\% & 62.6\% & 13.0\% & 36.5\%\\
 \hline
\end{tabular}
\end{center}
\caption{Historically reported relative change from the constrained MLE
when estimating \(\Omega\), with \(p=500\).  Positive values indicate
lower loss.  The procedures are the DAG completion of posterior-mean
precision coordinates, the full-Stein action \(\bar\Sigma^{-1}\), and
the transformed Cholesky-mode plug-in.  Saved outputs and Monte Carlo
uncertainty are unavailable.}
\label{tab:Estimation-p500-Precision-detail}
\end{table}

The historical results are highly sensitive to the hyperparameters and
to the generating edge probability.  In particular, Supplemental
Section~D reports severe deterioration at edge probability \(0.02\);
there is no uniform dominance over the MLE.  The supplement also
contains one Gaussian-contamination experiment.  Better performance
than the MLE in that single design is evidence about shrinkage under
that contamination, not a robustness theorem or a general guarantee
under misspecification.

\subsection{Historical molecular-network example}
The original analysis used the data of \cite{Sachs2003}, containing
measurements of \(p=11\) proteins and phospholipids on \(n=7466\)
cells, and imposed an ordering derived from the network reported in that
study.  That reported network is a scientific benchmark, not known
ground truth; accordingly, the terms true and false positive below mean
agreement or disagreement with that benchmark, not validated causal
truth.  The exact input data and preprocessing object used to produce
Figure~\ref{TrueData} are not present in the companion archive, so the
analysis is not independently reproduced here.

With \(\kappa=0.1\) for LassoDAG and \(b=3,c=1\) for DAG-W, the
historically reported rates are 78.95\% sensitivity and 52.78\%
specificity for LassoDAG, and 94.74\% sensitivity and 47.22\%
specificity for DAG-W.  Relative to a 19-edge benchmark among 55
ordered pairs, these rates correspond to approximately 15 benchmark
edges plus 17 additional edges for LassoDAG (precision \(46.9\%\)) and
18 benchmark edges plus 19 additional edges for DAG-W (precision
\(48.6\%\)).  Both selected graphs are therefore dense.  The result
shows a sensitivity--false-discovery tradeoff under the assumed
ordering; it does not by itself validate additional direct or indirect
molecular mechanisms.

\begin{figure}[htbp]
\centering
\SachsNetworkPanel{Benchmark DAG}{%
  PLCG/PIP3,PLCG/PIP2,PLCG/PKC,PKA/RAF,PKA/MEK,PKA/ERK,
  PKA/AKT,PKA/JNK,PKA/P38,PIP3/PIP2,PIP3/AKT,PIP2/PKC,
  PKC/RAF,PKC/MEK,PKC/JNK,PKC/P38,RAF/MEK,MEK/ERK,ERK/AKT}{}%
\hfill
\SachsNetworkPanel{LassoDAG estimate}{%
  PLCG/PIP3,PLCG/PIP2,PLCG/PKC,PKA/RAF,PKA/ERK,PKA/AKT,
  PKA/JNK,PIP3/PIP2,PIP3/AKT,PIP2/PKC,PKC/RAF,PKC/JNK,
  PKC/P38,RAF/MEK,ERK/AKT}{%
  PLCG/PKA,PLCG/RAF,PLCG/MEK,PLCG/ERK,PLCG/AKT,PLCG/JNK,
  PLCG/P38,PKA/PKC,PIP3/PKC,PIP3/ERK,PKC/ERK,PKC/AKT,
  MEK/AKT,MEK/P38,AKT/JNK,AKT/P38,JNK/P38}%
\hfill
\SachsNetworkPanel{DAG-W estimate}{%
  PLCG/PIP3,PLCG/PIP2,PLCG/PKC,PKA/RAF,PKA/MEK,PKA/ERK,
  PKA/AKT,PKA/JNK,PKA/P38,PIP3/PIP2,PIP3/AKT,PIP2/PKC,
  PKC/RAF,PKC/MEK,PKC/JNK,PKC/P38,RAF/MEK,ERK/AKT}{%
  PLCG/PKA,PLCG/RAF,PLCG/MEK,PLCG/ERK,PLCG/AKT,PLCG/JNK,
  PLCG/P38,PKA/PKC,PIP3/PKC,PIP3/MEK,PKC/ERK,PKC/AKT,
  RAF/AKT,MEK/AKT,MEK/P38,ERK/JNK,AKT/JNK,AKT/P38,JNK/P38}%
\caption{Historical selected graphs compared with the network benchmark
from \cite{Sachs2003}.  Blue edges agree with the benchmark and red
edges do not; colors do not establish causal truth.  The analysis could
not be reproduced from the public archive because the exact input
object is absent.}
\label{TrueData}
\end{figure}

Supplemental Section~D also preserves a historical call-center
prediction example based on data used in \cite{Levina08} and
\cite{rajaratnam08*}.  Its R script is available, but the external data
file and saved fitted model are not.  The numerical comparison is
therefore not independently verified in this revision.

%%%%%%%%%%%%%%%%%%%%%%%%%%%%%%%%%%%%5

\section{Closing remarks}\label{sec:closing}
The DAG--Wishart construction is valid for every DAG once an ordering
is fixed.  Its main tractable objects are the nodewise Cholesky
parameters: their independence yields the exact normalizing constant,
conjugate posterior, marginal likelihood, simulation scheme, and the
full posterior means in
\eqref{eq:full-precision-mean} and
\eqref{eq:cov-mean-offdiag}--\eqref{eq:cov-mean-diag}.  The induced
Lebesgue densities on incomplete coordinates are also explicit.

Three limitations are essential.  First, nonlinear completion does not
commute with expectation.  Second, a mode computed with respect to
Lebesgue measure in Cholesky coordinates is not automatically a
covariance- or precision-coordinate MAP.  Third, the model-selection
procedure in this paper assumes the ordering and therefore does not
resolve Markov-equivalence or unknown-order learning.

Subsequent work established high-dimensional graph-selection and
estimation consistency under explicit assumptions
\cite{CaoKhareGhosh2019}.  Other work showed that arbitrary
DAG--Wishart hyperparameters can assign different marginal likelihoods
to Markov-equivalent DAGs and constructed compatible specifications
\cite{PelusoConsonni2020}.  Those results should be consulted when the
inferential target extends beyond the fixed-order model class studied
here.

The numerical examples preserved in this preprint are historical
evidence.  A new empirical release should include immutable input data,
package and hardware manifests, all seeds, replicate-level outputs,
uncertainty intervals, run times, and comparisons that report false
discoveries as well as sensitivity.

%%%%%%%%%%%%%%%%%%%%

%%%%%%%%%%%%%%%%%%%%
\newpage
%%%%%%%%%%%%%%%%%%%%%%%%%%%% Suplementary Sections

\section*{Supplemental Section A: Graph theory and Gaussian DAG preliminaries}
\medskip
\subsection*{ Graph theoretic notation and terminology}\label{sec:appendixa}
A graph $\G$ is a pair of objects $(V, E)$, where $V$ and $E$ are two disjoint finite sets representing, respectively, the vertices and the edges of $\G$. Each edge in $E$ is either an ordered pair $(i,j)$ or an unordered pair $\{ i,j\}$, for some $i,j\in V$. An edge $(i,j)\in E$ is called directed where $i$ is said to be a parent of $j$, and $j$ is said to be a child of $i$, when $i\neq j$. We write this as $i\rightarrow j$. The set of parents of $i$ is denoted by $\mathrm{pa}(i)$, and the set of children of $i$ is denoted by $\mathrm{ch}(i)$. The family of $i$ is $\mathrm{fa}(i)=\mathrm{pa}(i)\cup \{i\}$. An edge $\{i,j\}\in E$ is called undirected where $i$ is said to be a neighbor of $j$, or $j$ a neighbor of $i$, when $i\neq j$. We write this $i\sim_{\G}j$. The set of all neighbors of $i$ is denoted by $\mathrm{ne}(i)$. We say $i$ and $j$ are adjacent if there exists  either a directed or an undirected edge between them. A loop in $\G$ is an ordered pair 
$(i,i)$, or an unordered pair $\{i,i\}$ in $E$. For ease of notation, we
include all loops in \(E\), although figures omit them.

\indent We say that the graph $\G'=(V', E')$ is a subgraph of  $\G=(V, E)$,  denoted by  $\G'\subset \G$,  if  $V'\subset V$  and  $ E'\subset  E$. In addition, if  $\G'\subset \G$  and  $ E'=(V'\times V')\cap E$, we say that $\G'$ is an {induced}  subgraph of  $\G$. We shall consider only induced subgraphs in what follows. For a subset  $A\subset V$, the induced subgraph  $\G_A=(A, (A\times A)\cap E)$ is said to be the graph {induced} by $A$. A graph  $\G$  is called {complete} if every pair of vertices are adjacent. A  clique  of  $\G$  is an induced complete subgraph  of  $\G$ that is not a subset of any other induced complete subgraph of  $\G$. More simply,  a subset  $A\subset V$  is  called a clique if the induced subgraph  $\G_A$  is a clique of $\G$. The set of the cliques of $\G$ is denoted by $\mathscr{C}_{\G}$.

\indent A path in $\G$ of length $n\geq1$ from $i$ to $j$ is a sequence
$i_0=i,\ldots,i_n=j$ such that each consecutive pair is adjacent; it is
simple when $i_0,\ldots,i_n$ are distinct.  In a directed graph, the path is
directed from $i$ to $j$ when $i_{\nu-1}\to i_\nu$ for every
$\nu=1,\ldots,n$.  We say that $i$ leads to $j$, written
$i\longmapsto j$, when such a directed path exists.  A graph is connected
when every distinct pair is joined by a path in its underlying undirected
graph.  A cycle of length $n\geq2$ is a sequence
$i_0,i_1,\ldots,i_{n-1},i_n=i_0$ in which
$i_0,\ldots,i_{n-1}$ are distinct and consecutive vertices are adjacent; it
is directed when every edge is oriented $i_{\nu-1}\to i_\nu$.

\indent An undirected graph $\G=(V,\mathscr E)$ is decomposable (or
chordal) if it has no induced cycle of length at least four.  A directed
graph $\D=(V,E)$ is acyclic when it has no directed cycle apart from the
notational loops.  Its skeleton $\D^u$ is obtained by replacing directed
edges by undirected edges.  An immorality is an induced configuration
$i\to k\leftarrow j$ with nonadjacent parents $i$ and $j$.  The moral graph
$\D^{\mathrm m}$ is obtained by joining every pair of co-parents and then
dropping all arrowheads.  A DAG is perfect when it has no immoralities,
equivalently when the parents of every vertex are pairwise adjacent.  Every
decomposable undirected graph admits a perfect acyclic orientation
\cite{lauritzen96*}.

\indent Given a DAG, the set of ancestors of a vertex $j$, denoted by $\mathrm{an}(j)$, is the set of those vertices $i$  such that $i\longmapsto j$. Similarly, the set of descendants of a vertex $i$, denoted by  $\mathrm{de}(i)$, is the set of  those  vertices  $j$  such that  $i\longmapsto j$.  The   set of non-descendants of  $i$  is  $\mathrm{nd}(i)=V\setminus\left(\mathrm{de}(i)\cup \{i\}\right)$. A set  $A\subseteq V$  is called ancestral when $A$  contains the parents of its members. The smallest ancestral set containing the subset $B$ of $V$ is denoted by $\mathrm{An}(B)$. 

\subsection*{Markov properties for DAG models}\label{sub:mpfdag}
Let $V$ be a finite set of indices and $(X_i)_{i\in V}$ a collection of random variables, where each $X_i$ is a random variable on the probability space $\mathcal{X}_i$. Let the probability space $\mathcal{X}$ be defined as the product space $\mathcal{X}=\times_{i\in V}\mathcal{X}_i$. Now let $\D=(V, E)$ be a DAG.  For simplicity, and without loss of generality, we always assume that the given DAG $\D$ is connected and the edge set $E$  contains all the loops $(i,i), i\in V$.  We say that a probability distribution  $P$  on  $\mathcal{X}$  has the recursive factorization property w.r.t. $\D$, denoted by  DF (the directed factorization property), if  there are  $\sigma$-finite  measures  $\mu_{i}$  on  $\mathcal{X}_i$  and non-negative functions   $k^i(x_i, x_{pa(i)})$, referred to as {kernels}, defined on   $\mathcal{X}_{fa(i)}$ such that
\[
\int k^i(y_i, x_{pa(i)})d\mu_i(y_i)=1,\quad\forall i\in V,
\]
and  $P$  has a density  $p$,  w.r.t.  the product measure  $\mu=\otimes_{i\in V}\mu_i$, given by
\[
p(x)=\prod_{i\in V}k^i(x_i,x_{pa(i)}).
\]
In this case, each kernel  $k^i(x_i,x_{pa(i)})$  is in fact a version of  $p(x_i|x_{pa(i)})$, the conditional distribution of $X_i$  given $X_{pa(i)}$.
An immediate consequence of this definition is the following lemma.

\begin{lemma}\cite{lauritzen96*}\label{dfimoralf}
If  $P$ admits a recursive factorization w.r.t. the directed graph $\D$, then it also admits a factorization w.r.t. the undirected graph  $\D^{\mathrm{m}}$,  and,  consequently, obeys the global Markov property\footnote{see \cite{lauritzen96*} for definition.} w.r.t. $\D^{\mathrm{m}}$.
\end{lemma}
\begin{proof}
Note that for each vertex $i\in V$ the set $fa(i)$ is a complete subset of $\D^{\mathrm{m}}$. Thus if we define  $\psi_{fa(i)}(x_{fa(i)})=k^{i}(x_{i}, x_{pa(i)})$, then  $p(x)= \prod_{i\in V}p(x_i|x_{pa(i)}) = \prod_{i\in V}k^{i}(x_i,x_{pa(i)}) = \prod_{i\in V}\psi_{fa(i)}(x_{fa(i)})$. Therefore, $P$ admits a factorization w.r.t. $\D^{\mathrm{m}}$ and by proposition 3.8 in \cite{lauritzen96*} it also obeys the global Markov property w.r.t. $\D^{\mathrm{m}}$.
\end{proof}

\indent
Another direct implication of the DF property is that if  $P$  admits a recursive factorization w.r.t.  $\D$, then, for each ancestral set $A$, the marginal distribution   $P_A$   admits a recursive factorization w.r.t.  the induced graph  $\D_A$. Combining this result with Lemma \ref{dfimoralf} we obtain the following: if \(P\) admits a recursive factorization w.r.t. \(\D\), then $A\Perp B| S\  [P]$ whenever  $A$  and  $B$ are separated by  $S$  in  $(\D_{An(A\cup B\cup S)})^{\mathrm{m}}$. We call this property the directed global Markov property, DG, and any distribution that satisfies this property is said to be a directed Markov field over $\D$. For DAGs the directed Markov property plays the same role as the global Markov property does for undirected graphs: it recovers the conditional independence relations encoded by the directed graph.

\indent We now introduce below another Markov property for DAGs. A  distribution   $P$  on  $\mathcal{X}$ is said to obey the {directed local Markov property} (DL) w.r.t. $\D$  if  for each  $i\in V$
  \[
  i\Perp \bigl(nd(i)\setminus pa(i)\bigr)\mid pa(i).
  \]
  For a given DAG \(\D\), a parent ordering relabels the vertices as
  \(1,\ldots,|V|\) so that
  \(pa(i)\subseteq\{i+1,\ldots,|V|\}\) for every \(i\).  Every DAG has
  at least one such ordering, although it need not be unique.  We say
  that \(P\) obeys the parent-ordered Markov property (PO) with respect
  to \(\D\) if, for every vertex \(i\),
 \[
  i\Perp \left\{i+1,\ldots, |V|\right\}\setminus pa(i)|pa(i).
  \]

\indent If \(P\) has a density with respect to \(\mu\), then the four
directed Markov properties DF, DG, DL, and PO are equivalent under the
usual regularity conditions \cite{lauritzen96*}.

\subsection*{Linear recursive properties of Gaussian DAGs}\label{subs: Lrp}
Let $\mathbf{x}=\left(x_1,\ldots,x_{p}\right)^{\top}$  be  a random vector in  $\R^{p}$  with    the multivariate distribution  $\rm{N}_{p}(0, \Sigma)$.  Consider the system of linear recursive regression equations:
\[
\scalebox{.85}{$\displaystyle
\begin{array}{rclcrcl}
x_1-\beta_{12}x_2-\beta_{13}x_3-\cdots-\beta_{1p}x_{p}&=&\epsilon_1 
&\text {or equivalently} &
x_1&=&\beta_{12}x_2+\beta_{13}x_3+\cdots+\beta_{1p}x_{p}+\epsilon_1\\
x_2-\beta_{23}x_3-\cdots-\beta_{2p}x_{p}&=&\epsilon_2
&\quad&
x_2&=&\beta_{23}x_3+\cdots+\beta_{2p}x_{p}+\epsilon_2\\
&\vdots& &&&\vdots& \\
x_{p}&=&\epsilon_{p}
&\quad&
x_{p}&=&\epsilon_{p},
\end{array}$}
\]

where \(\beta_{ij}\) is the regression coefficient of \(x_j\)
(\(j>i\)) in the regression of \(x_i\) on its predecessors.  It is
zero when \(j\notin\mathrm{pa}(i)\).  The residuals \(\epsilon_i\) are
mutually independent centered normal variables with variance
\(\sigma^2_{ii\mid pa(i)}\).  The equations can be written
\(B\mathbf{x}=\boldsymbol{\epsilon}\), where \(B\) is upper triangular:
\[
B=
\left(
\begin{matrix}
1& -\beta_{12}&\ldots& -\beta_{1p}\\
  0&1 &\ldots &-\beta_{2p}\\
  0&\ldots&\ddots&\vdots \\
 0 &\ldots&0 & 1\\
  \end{matrix}
\right)
,\quad \mathbf{x}=\left(\begin{matrix}x_1\\x_2\\\vdots\\ x_{p}\end{matrix}\right)
\quad\text{and}\quad\boldsymbol{\epsilon}=\left(\begin{matrix}\epsilon_1\\\epsilon_2\\\vdots\\ \epsilon_{p}\end{matrix}\right).
\]
From this we obtain:
\begin{align}
\notag
&\mathbb{V}ar[B\mathbf{x}]=\mathbb{V}ar[\mathbf{\epsilon}]\\
\notag\Rightarrow& B\Sigma B^{\top}=diag(\sigma^2_{1|pa(1)}, \ldots, \sigma^2_{p-1|pa(p-1)}, \sigma^2_{pp})=: D\\
\notag
\Rightarrow& \Sigma= B^{-1}D(B^{\top})^{-1}\\
\Rightarrow& \Sigma^{-1}= B^{\top}D^{-1}B.\label{eq: choleskyd}
\end{align}
Thus \(L=B^\top\) gives
\(\Sigma^{-1}=LD^{-1}L^\top\), with
\(L_{\mathrm{pa}(i),i}=-\beta_i\), exactly as specified in
Section~\ref{subs:d_parameterization}.  The Gaussian distribution is
Markov with respect to \(\D\) precisely when the corresponding
non-parent entries of \(L\) vanish \cite{wermuth80}.
%

%%%%%%%%%%%%%%%%%%%%%%%%%%%%%%%%%%%%%%%

\section*{Supplemental Section B: Properties of the class of DAG-Wishart distributions}

\subsection{Deriving the closed form expression for the DAG-Wishart $\pi_{U,\alpha}^{\Theta_{\D}}$}

\begin{theorem} \label{nrmcnstdml}
Let $dL := \prod_{(i,j) \in E, i>j} dL_{ij}$ and $dD := \prod_{i=1}^p dD_{ii}$ denote, respectively,  the canonical Lebesgue measures on $\mathcal{L}_{\D}$  and  $\R^p_{+}$ and let $pa_i:=|pa(i)|$. Then,
$$
 \int_{\Theta_{\D}}\exp\{- \frac{1}{2}\tr( LD^{-1}L^t U)\}
\prod_{i=1}^p D^{-\frac{1}{2}\alpha_i}_{ii} dL dD < \infty
$$
if and only if
$$\alpha_i > pa_i + 2 \;\;\; \forall i = 1,\ldots, p. $$  

\noindent Furthermore, in this case 
\begin{equation} \label{nrmlcnsten}
z_{\D} (U, \alpha) = \prod_{i=1}^p \frac{\Gamma \left( \frac{\alpha_i}{2} -
\frac{pa_i}{2} - 1 \right) 2^{\frac{\alpha_i}{2} - 1}
(\sqrt{\pi})^{pa_i} \det( U_{\prec i\succ})^{\frac{\alpha_i}{2} -
\frac{pa_i}{2} - \frac{3}{2}}}{\det( U_{\preceq i\succeq})^{\frac{\alpha_i}{2} - \frac{pa_i}{2} - 1}}.
\end{equation}
\end{theorem}
\begin{proof}Let us first simplify the expression by integrating out the terms involving  $D_{ii}$'s.
\begin{eqnarray*}
& & \int \exp\left\{- \frac{1}{2}\tr \left(\left(LD^{-1}L^{\top}\right) U \right) \right\}\prod^{p}_{i=1}D_{ii}^{-\frac{1}{2}\alpha_i} dL dD\\
&=& \int \exp\left\{-\frac{1}{2} \tr \left( D^{-1}\left( L^{\top}U L\right) \right)\right\}
\prod_{i=1}^{p}  D^{-\frac{1}{2}\alpha_i}_{ii} dL dD\\
&=& \int \exp\left\{- \frac{1}{2}\sum^{p}_{i=1}D^{-1}_{ii}( L^{\top} U L)_{ii}\right\}
\prod^{p}_{i=1}D_{ii}^{-\frac{1}{2}\alpha_i} dD dL\\
&=& \int\left( \prod_{i=1}^{p} \int\exp\left\{- \frac{1}{2}D^{-1}_{ii}( L^{\top} U L)_{ii}\right\}
D_{ii}^{-\frac{1}{2}\alpha_i} dD_{ii}\right) dL\\
&=& \int \prod_{i=1}^{p} \frac{\Gamma \left( \frac{\alpha_i}{2} - 1 \right)
2^{\frac{\alpha_i}{2} - 1}}{\left( \left( L^{\top} U L \right)_{ii}
\right)^{\frac{\alpha_i}{2} - 1}} dL \qquad(\text{ iff  } \alpha_i > 2 \; \forall \; i =
1, 2, \cdots, p)\\
&=& \int \prod_{i=1}^{p} \frac{\Gamma \left( \frac{\alpha_i}{2} - 1 \right)
2^{\frac{\alpha_i}{2} - 1}}{\left( (L_{\cdot i })^{\top} U L_{\cdot i}
\right)^{\frac{\alpha_i}{2} - 1}} dL\\
&=& \int \prod_{i=1}^{p} \frac{\Gamma \left( \frac{\alpha_i}{2} - 1 \right)
2^{\frac{\alpha_i}{2} - 1}}{\left(\left(\begin{matrix}1& L^{\top}_{\prec i]}\end{matrix}\right)\left(\begin{matrix}U_{ii}& U_{[i\succ}\\U_{\prec i]}& U_{\prec i\succ}\end{matrix}\right) \left(\begin{matrix} 1\\ L_{\prec i]} \end{matrix}\right)
\right)^{\frac{\alpha_i}{2} - 1}} dL\\
&=& \prod_{i=1}^{p} \int_{\R^{pa_i}} \frac{\Gamma \left( \frac{\alpha_i}{2} - 1 \right)
2^{\frac{\alpha_i}{2} - 1}}{\left(\left(\begin{matrix}1& L^{\top}_{\prec i]}\end{matrix}\right)\left(\begin{matrix}U_{ii}& U_{[i\succ}\\U_{\prec i]}& U_{\prec i\succ}\end{matrix}\right) \left(\begin{matrix} 1\\ L_{\prec i]} \end{matrix}\right)
\right)^{\frac{\alpha_i}{2} - 1}} dL_{\prec i]}.\qquad \text{eqn(A)}
\end{eqnarray*}
We now show how in general one can evaluate an integral of the form
\[
\int_{\R^d} \frac{d\mathbf{x}}{\left(\left(\begin{matrix}1& \mathbf{x}^{\top} \end{matrix}\right)\left(\begin{matrix} a  & \mathbf{b}^{\top}\\  \mathbf{b}&A\end{matrix}\right)\left(\begin{matrix} 1\\ \mathbf{x} \end{matrix}\right)\right)^{\gamma}},
\]
where the block matrix, formed by \(a\in\R\),
\(\mathbf{b}\in\R^d\), and a \(d\times d\) matrix \(A\), is positive
definite.
\noindent
To simplify the integral, proceed in two steps.

\noindent \(1)\) By \cite[page 16]{diaconis2008},
\[
\int_{\mathbb{R}} \frac{1}{(1 + x^2)^\gamma} dx = \begin{cases}
\frac{\sqrt{\pi} \Gamma \left( \gamma - \frac{1}{2} \right)}{\Gamma(\gamma)} &
\gamma > \frac{1}{2}, \cr
\infty & \mbox{otherwise}.
\end{cases}
\]
\noindent
By repeated application, we can generalize the above formula to
\[
\int_{\mathbb{R}^d} \frac{1}{({\bf x}^{\top} {\bf x} + 1)^\gamma} d {\bf x} = \begin{cases}
\frac{(\sqrt{\pi})^d \Gamma \left( \gamma - \frac{d}{2} \right)}{\Gamma(\gamma)} & \gamma
> \frac{d}{2}, \cr
\infty & \mbox{otherwise}.
\end{cases}
\]
\noindent \(2)\) Now consider the general integral
\[
	\int_{\R^d} \dfrac{d\mathbf{x}}{\left(\left(\begin{matrix}1& \mathbf{x}^{\top} \end{matrix}\right)\left(\begin{matrix} a  & \mathbf{b}^{\top}\\  \mathbf{b}&A\end{matrix}\right) \left(\begin{matrix} 1\\ \mathbf{x} \end{matrix}\right)\right)^{\gamma}} .
\]
\noindent
Making the linear transformation ${\bf y} = A^{\frac{1}{2}} {\bf x}
+ A^{-\frac{1}{2}}\mathbf{b}$ it follows that for $\gamma > \frac{d}{2}$,
\begin{eqnarray}
\int_{\R^d} \dfrac{d\mathbf{x}}{\left(\left(\begin{matrix}1& \mathbf{x}^{\top} \end{matrix}\right)\left(\begin{matrix} a  & \mathbf{b}^{\top}\\  \mathbf{b}&A\end{matrix}\right)\left(\begin{matrix} 1\\ \mathbf{x} \end{matrix}\right)\right)^{\gamma}}
&=& \frac{1}{\det(A)^{\frac{1}{2}}} \int_{\mathbb{R}^d} \dfrac{1}{\left( {\bf y}^{\top} {\bf y} +  a - \mathbf{b}^{\top} A^{-1}
\mathbf{b} \right)^\gamma} d {\bf y} \nonumber\\
&=& \dfrac{(\sqrt{\pi})^d \Gamma \left( \gamma - \frac{d}{2} \right)}{\Gamma(\gamma)
\det(A)^{\frac{1}{2}} \left( a - \mathbf{b}^{\top} A^{-1} \mathbf{b} \right)^{\gamma - \frac{d}{2}}}.
\label{intgltrnfm}
\end{eqnarray}
\noindent
Applying the result from  \eqref{intgltrnfm} to the desired integral in (A) we obtain
\begin{eqnarray*}
z_{\D} (U, {\bf \alpha}) &=&\prod_{i=1}^{p} \int_{\R^{pa_i}} \frac{\Gamma \left( \frac{\alpha_i}{2} - 1 \right)
2^{\frac{\alpha_i}{2}-1}}{\left(\left(\begin{matrix}1& L^{\top}_{\prec i]}\end{matrix}\right)\left(\begin{matrix}U_{ii}& U_{[i\succ}\\U_{\prec i]}& U_{\prec i\succ}\end{matrix}\right) \left(\begin{matrix} 1\\ L_{\prec i]} \end{matrix}\right)
\right)^{\frac{\alpha_i}{2} - 1}} dL_{\prec i]}\\
&=& \prod_{i=1}^{p} \frac{\Gamma \left( \frac{\alpha_i}{2} -
\frac{pa_i}{2} - 1 \right) 2^{\frac{\alpha_i}{2} - 1}
(\sqrt{\pi})^{pa_i} \det(U_{\prec i\succ})^{\frac{\alpha_i}{2} -
\frac{pa_i}{2} - \frac{3}{2}}}{\det(U_{\preceq i\succeq})^{\frac{\alpha_i}{2} -
\frac{pa_i}{2} - 1}},
\end{eqnarray*}
\noindent
where $\det(U_{\prec i \succ}) = 1$ whenever $pa(i) = \emptyset$. It is easily seen that $z_{\D} (U, {\bf \alpha})$ is finite  iff  $\alpha_i> pa_i+2 $  for each  $i=1,\ldots,p.$
%\noindent
\end{proof}

%%%%%%%%%%%%%%%%%%%%%%%%%%%%%%%%%%%%%%%%%%
\subsection*{The hyper-Markov properties of the DAG-Wishart}\label{sub:hydag}
\begin{theorem}\label{hydag}
Let $\D$ be an arbitrary DAG and $(D,L) \sim \pi^{\Theta_{\D}}_{U, {\bf \alpha}}$. Then $\{(D_{ii}, L_{\prec i]}): i=1,\ldots, p\}$ are mutually independent. Moreover,
\begin{equation}\label{dii}
D_{ii}\sim IG(\frac{\alpha_i}{2}-\frac{pa_i}{2}-1,  \frac{1}{2}U_{ii|\prec i\succ}) ,~\text{and}
\end{equation}
\begin{equation}\label{ind}
L_{\prec i]}|D_{ii}\sim \mathrm{N}_{pa_i}(-U^{-1}_{\prec i\succ}U_{\prec i]}, D_{ii}U^{-1}_{\prec i\succ} ).
\end{equation}
\end{theorem}

\begin{proof}
The coordinate map
\begin{equation}\label{oto}
\phi:(D,L)\longmapsto
\mathop{\times}_{i\in V}(D_{ii},\ell_i),
\qquad \ell_i=L_{\prec i]},
\end{equation}
has unit Jacobian.  It is a vectorization of the Cholesky parameters,
not of the regression coefficients; the latter are
\(\beta_i=-\ell_i\).  Cyclicity of trace and completion of the square
give
\begin{align*}
&\tr((LD^{-1}L^{\top})U)= \tr(D^{-1}L^{\top}UL) =\sum_{i\in V}D^{-1}_{ii}(L^{\top}UL)_{ii}\\
&=\sum_{i\in V}D^{-1}_{ii}(\sum_{k,l\in V}L_{ki}U_{kl}L_{li})\\
&=\sum_{i\in V}D^{-1}_{ii}\left(\begin{matrix}1\\ L_{\prec i]}\end{matrix}\right)^{\top}\left(\begin{matrix}U_{ii}& U_{[i\succ}\\U_{\prec i]}& U_{\prec i\succ} \end{matrix}\right)\left(\begin{matrix}1\\L_{\prec i]}  \end{matrix}\right)\\
&=\sum_{i\in V}D^{-1}_{ii}(U_{ii}+ L_{\prec i]}^{\top}U_{\prec i]}+U_{[i\succ}L_{\prec i]}+L^{\top}_{\prec i]}U_{\prec i\succ}L_{\prec i]})\\
&=\sum_{i\in V}\left(D^{-1}_{ii}(L_{\prec i]}+U^{-1}_{\prec i\succ}U_{\prec i]})^{\top} U_{\prec i\succ} (L_{\prec i]}+U^{-1}_{\prec i \succ}U_{\prec i]}) +  D^{-1}_{ii}U_{ii|\prec i\succ}\right).
\end{align*}
Therefore, the density in the coordinates
\((\lambda_i,\ell_i)=(D_{ii},L_{\prec i]})\) is
\begin{equation}\label{D-parameterp}
z_{\D}(U,\alpha)^{-1}
\prod_{i\in V}
\left[
\lambda_i^{-\alpha_i/2}
\exp\left\{-\frac{1}{2\lambda_i}
\left(
U_{ii\mid\prec i\succ}
+(\ell_i+U_{\prec i\succ}^{-1}U_{\prec i]})^\top
U_{\prec i\succ}
(\ell_i+U_{\prec i\succ}^{-1}U_{\prec i]})
\right)\right\}
\right].
\end{equation}
The product proves mutual independence across vertices.  The \(i\)th
factor is the product of the inverse-gamma density in \eqref{dii} and
the conditional Gaussian density in \eqref{ind}, including their
normalizing constants.  This proves both stated laws.
\end{proof}

%%%%%%%%%%%%%%%%%%%%%%%%%%%%%%%%%%%%%%%%%%%%%%%%%%%%%%%%%%%%%%%%%%%%%%%%%%%%%%%%%%%%%%%%%%%%%%
\begin{corollary}\label{corr: tdisn}
Let \((L,D)\sim\pi^{\Theta_{\D}}_{U,\alpha}\) and set
\(\nu_i=\alpha_i-pa_i-2>0\).  The columns \(L_{\prec i]}\) are mutually
independent and
\[
L_{\prec i]}\sim
t_{\nu_i}\!\left(
-U_{\prec i\succ}^{-1}U_{\prec i]},
\frac{U_{ii\mid\prec i\succ}}{\nu_i}
U_{\prec i\succ}^{-1}\right).
\]
Equivalently, their joint density is
\[
\prod_{i=1}^{p}c_i
\left[
\frac12\left\{
U_{ii\mid\prec i\succ}
+(L_{\prec i]}+U_{\prec i\succ}^{-1}U_{\prec i]})^\top
U_{\prec i\succ}
(L_{\prec i]}+U_{\prec i\succ}^{-1}U_{\prec i]})
\right\}
\right]^{-\alpha_i/2+1},
\]
where
\begin{equation}\label{eq:c_i}
c_i=
\frac{\det(U_{\prec i\succ})^{1/2}
(U_{ii\mid\prec i\succ})^{\alpha_i/2-pa_i/2-1}
\Gamma(\alpha_i/2-1)}
{2^{\alpha_i/2-1}\pi^{pa_i/2}
\Gamma(\alpha_i/2-pa_i/2-1)}.
\end{equation}
\end{corollary}

\begin{proof}
Multiply the conditional Gaussian density in \eqref{ind} by the
inverse-gamma density in \eqref{dii} and integrate \(D_{ii}\) over
\((0,\infty)\).  The integral is
\(\Gamma(\alpha_i/2-1)u_i^{-\alpha_i/2+1}\), where
\[
u_i=\frac12\left\{
U_{ii\mid\prec i\succ}
+(L_{\prec i]}+U_{\prec i\succ}^{-1}U_{\prec i]})^\top
U_{\prec i\succ}
(L_{\prec i]}+U_{\prec i\succ}^{-1}U_{\prec i]})
\right\}.
\]
Collecting constants gives \eqref{eq:c_i}; comparison with the
multivariate \(t\) density gives the stated degrees of freedom and
scale.
\end{proof}
%%%%%%%%%%%%%%%%%%%%%%%%%%%%%%%%%%%%%%%%%%%%%%%%%%%%%%%%%%%%%%%%%%%%%%%%%%%%%%%%%%%%%
%%%%%%%%%%%
\subsection*{The Posterior distribution of the DAG-Wishart}

\begin{proposition}\label{lem:pd}
Let $\D$ be an arbitrary DAG and let ${\bf Y}_1, {\bf Y}_2, \cdots, {\bf Y}_n$ be an i.i.d. sample from $\rm{N}_{p} ({\bf 0}, (L^{-1})^{\top}DL^{-1})$, where
$(D, L) \in \Theta_{\D}$. Let $S = \frac{1}{n} \sum_{i=1}^n {\bf Y}_i {\bf Y}_i^{\top}$ denote the empirical
covariance matrix. If the prior distribution on $(D, L)$ is $\pi_{U,\alpha}^{\Theta_{\D}}$, then the posterior distribution
of $(D, L)$ is given by $\pi_{\widetilde{U}, \widetilde{\bf \alpha}}^{\Theta_{\D}}$, where $\widetilde{U} = nS + U$ and $\widetilde{\bf \alpha} = (n + \alpha_1, n + \alpha_2, \cdots, n + \alpha_{p})$.
\end{proposition}
\begin{proof}
The likelihood of the data is given as follows:
\[
f({\bf y}_1, {\bf y}_2, \cdots, {\bf y}_n \mid L, D) = \frac{1}{(\sqrt{2 \pi})^{np}}\exp\left\{-
\frac{1}{2}\tr \left(LD^{-1}L^{\top} (nS) \right)\right\}\det(D)^{-\frac{1}{2}n}.
\]
\noindent When using $\pi_{U,\alpha}^{\Theta_{\D}}$ as the prior for $(D, L)$, the posterior distribution of $(D, L)$ given the
data $({\bf Y}_1, {\bf Y}_2, \cdots, {\bf Y}_n)$ is given by
\begin{equation}\label{eq:wishart_on_theta}
\begin{aligned}
\pi_{U,\alpha}^{\Theta_{\D}}(L,D\mid{\bf Y}_{1:n})
&\propto
\exp\left\{-\frac12\tr\!\left(LD^{-1}L^\top(nS+U)\right)\right\}\\
&\quad{}\times\prod_{i=1}^{p}D_{ii}^{-(n+\alpha_i)/2},
\qquad (D,L)\in\Theta_{\D}.
\end{aligned}
\end{equation}

\noindent
Hence the functional form of the posterior density is the same as that of the prior density, i.e.,
$$
\pi_{U,\alpha}^{\Theta_{\D}} ( \cdot \mid {\bf Y}_1, {\bf Y}_2, \cdots, {\bf Y}_n) =
\pi_{\widetilde{U}, \widetilde{\bf \alpha}} ( \cdot ),
$$

\noindent
where $\widetilde{U} = nS + U$ and $\widetilde{\alpha}= (\alpha_1+n,\ldots,\alpha_{p}+n)$.
\end{proof}
\begin{Rem}\label{rem:unknown-mean}
If \(Y_r\stackrel{\mathrm{iid}}{\sim}\mathrm{N}_p(\mu,\Sigma)\) and
the prior for \(\mu\) is the independent flat prior
\(\pi(\mu)\propto1\), integrate \(\mu\) out and define
\[
\bar Y=n^{-1}\sum_{r=1}^nY_r,\qquad
T_c=\sum_{r=1}^n(Y_r-\bar Y)(Y_r-\bar Y)^\top .
\]
Then \(T_c\sim W_p(n-1,\Sigma)\), independently of \(\bar Y\), and the
posterior Cholesky law is
\[
\pi_{U+T_c,\,\alpha+(n-1)\mathbf{1}}^{\Theta_{\D}}.
\]
Equivalently, if \(S_c=T_c/(n-1)\), then
\((n-1)S_c\sim W_p(n-1,\Sigma)\).  This update depends on the stated
flat prior; a different prior or conditioning treatment for \(\mu\)
must be derived separately.
\end{Rem}
%%%%%%%%%%%%%%%%%%%%%%%%%%%%%%%%%%%%%%%%%%%%%%%%%%%%%%%%%%%%%%%%%%%

%%%%%%%%%%%%%%%%%%%%%%%%%%%%%%%%%%%%%%%%%%%%%%%%%%%%%%%%%%%%%%%%%%%%%%%%%%%%%%%%%%%%%
%%%%%%%%%%%%%%%%%%%%%%%%%%%%%%%%%%%%%%%%%%%%%%%%%%%%%%%%%%%%%%%%%%%%%%%%%%%%%%%%%%%%%
%%%%%%%%%%%%%%%%%%%%%%%%%%%%%%%%%%%%%%%%%%%%%%%%%%%%%%%%%%%%%%%%%%%%%%%%%%%%%%%%%%%%%
%%%%%%%%%%%%%%%%%%%%%%%%%%%%%%
\subsection*{ The Laplace transform of the DAG-Wishart} We start with computing the Laplace transform of  $\pi^{\Xi_{\D}}_{U,\alpha}$ by exploiting the results established in Theorem \ref{hydag}. First a preliminary result on the Laplace transform of a Gaussian inverse Gamma distribution is required.
\begin{lemma}\label{lem:laplace_t_normal-inv-gamma}
Suppose $(\lambda, \mathbf{x})$ is a random variable with Gaussian-inverse gamma distribution:
\begin{align*}
\mathbf{x}|\lambda&\sim \mathrm{N}_{p}(\mu, \lambda \Psi) ,\quad \mu\in \R^{p}, \Psi\in \mathrm{PD}_{p}(\R);\\
\lambda&\sim \mathrm{IG}(\nu, \eta).
\end{align*}
For \(u\in\R^p\), the Laplace transform
\(E[\exp\{-\xi\lambda-u^\top\mathbf{x}\}]\) is finite exactly when
\(\xi\geq\tfrac12u^\top\Psi u\).  When the inequality is strict, it equals
\[
\frac{2}{\Gamma(\nu)}\exp\{-u^{\top}\mu\}\left(\eta(\xi-\frac{1}{2}u^{\top}\Psi u)\right)^{\frac{1}{2}\nu} K_{\nu}\left(2\sqrt{\eta(\xi-\frac{1}{2}u^{\top}\Psi u)}\right),
\]
where \(K_{\nu}\) is the modified Bessel function of the second kind.
At equality, the continuous limit is \(\exp(-u^\top\mu)\).
\end{lemma}
\begin{proof} 
By definition, the Laplace transform of $(\lambda, \mathbf{x})$  at  $(\xi, u)\in \R\times \R^{p}$  is 
\begin{align*}\label{eq:Laplace_t_gaussian}
\notag
&\int\exp\{ -(\lambda\xi+ u^{\top}x)\}d\mathrm{N}_{p}(\mu, \lambda\Psi)(x)d\mathrm{IG}(\nu, \eta)(\lambda)\\
=&\int\exp\{ -\lambda\xi\}\left(\int\exp\{- u^{\top}x\}d\mathrm{N}_{p}(\mu, \lambda\Psi)(x)\right)d\mathrm {IG}(\nu, \eta)(\lambda)\\
\notag=&\int\exp\{ -\lambda\xi\}\exp\{-u^{\top}\mu+\frac{1}{2}\lambda u^{\top}\Psi u  \}d\mathrm{IG}(\nu, \eta)(\lambda)\\
=&\int\exp\{ -\lambda\xi\}\exp\{-u^{\top}\mu+\frac{1}{2}\lambda u^{\top}\Psi u \}\left( \frac{\eta^{\nu}}{\Gamma(\nu)}\exp\{ -\eta\lambda^{-1}\}\lambda^{-\nu-1}\right) d\lambda\\
=&\frac{\eta^{\nu}}{\Gamma(\nu)}\exp\{- u^{\top}\mu\}\int \exp\{-(\xi-\frac{1}{2}u^{\top}\Psi u)\lambda -\eta\lambda^{-1}\} \lambda^{-\nu-1}d\lambda\\
=&\frac{2\eta^{\nu}}{\Gamma(\nu)}\exp\{- u^{\top}\mu\}\left(\frac{\xi-\frac{1}{2}u^{\top}\Psi u}{\eta}\right)^{\frac{1}{2}\nu} K_{\nu}\left(2\sqrt{\eta(\xi-\frac{1}{2}u^{\top}\Psi u)}\right) \\
=&\frac{2}{\Gamma(\nu)}\exp\{- u^{\top}\mu\}\left(\eta(\xi-\frac{1}{2}u^{\top}\Psi u)\right)^{\frac{1}{2}\nu} K_{\nu}\left(2\sqrt{\eta(\xi-\frac{1}{2}u^{\top}\Psi u)}\right).
\end{align*}
 Note that  in computing the integral above we have used the fact that the Laplace transform of $
\mathrm{N}_{p}(\mu, \lambda\Psi)$ at  $u$ is  equal to  $\exp\{ -u^{\top}\mu+\frac{1}{2}\lambda u^{\top}\Psi u\}$. For computing the integral w.r.t. $d\lambda$  we use the Equation (9.42) in \cite[page 235]{temme97}.

\end{proof}
%------------------------------------
\medskip

\begin{proposition}
At test arguments \(\times_{i=1}^{p}(\xi_i,z_{\prec i]})\) satisfying
\(\xi_i>\tfrac12z_{\prec i]}^\top\Psi_{\prec i\succ}z_{\prec i]}\),
the Laplace transform of $\pi^{\Xi_{\D}}_{U,\alpha}$ is
\begin{eqnarray}
\mathcal{L}_{\Xi_{\D}}(\times_{i=1}^{p}(\xi_i, z_{\prec i]})):=2^{p}\prod_{i=1}^{p}   \frac{1}{\Gamma(r_{i})}\exp\{- z_{\prec i]}^{\top}\mu_{\prec i]}\}\left(\eta_{i}(\xi_{i}-\frac{1}{2}z_{\prec i]}^{\top}\Psi_{\prec i\succ} z_{\prec i]})\right)^{\frac{1}{2}r_{i}}\\
\times K_{r_{i}}\left(2\sqrt{\eta_{i}(\xi_{i}-\frac{1}{2}z_{\prec i]}^{\top}\Psi_{\prec i \succ} z_{\prec i]})}\right),
\end{eqnarray}
where $r_{i}=\frac{\alpha_{i}}{2}-\frac{pa_{i}}{2}-1$, \:
$\eta_{i}=\frac{1}{2}U_{ii| \prec i\succ}$, \:
$\mu_{\prec i]}=-U_{\prec i\succ}^{-1}U_{\prec i]}$, and \:
$\Psi_{\prec i \succ}=U_{\prec i\succ}^{-1}$.  Boundary arguments with
equality are obtained by taking the corresponding continuous limits.
\end{proposition}
\begin{proof}
Let
\(\times_{i=1}^{p}(\lambda_i,\ell_i)\sim
\pi^{\Xi_{\D}}_{U,\alpha}\).  Theorem~\ref{hydag} makes these
Gaussian--inverse-gamma pairs independent, so the result is the product
of Lemma~\ref{lem:laplace_t_normal-inv-gamma} applied with the displayed
parameters.
\end{proof}
%%%%%%%%%%%%%%%%%%%%%%%%%%%%%%%%%%%%%%%%%%%%%%%%%%%%%%
%%%%%%%%%%%%%%%%%%%%%%%%%%%%%%%%%%%%%%%%%%%%%%%%%%%%%%

Because the coordinate map \(\phi\) in \eqref{oto} has unit Jacobian,
the preceding proposition is also the joint Laplace transform of the
\emph{free} Cholesky coordinates
\(\{D_{ii},L_{\prec i]}:i=1,\ldots,p\}\).  The fixed unit diagonal of
\(L\) is not included as a random coordinate.  A matrix-trace version
can be obtained by setting
\(\xi_i=\Lambda_{ii}\) and taking \(z_{\prec i]}\) from the free
entries of the test matrix; any deterministic diagonal contribution
must be included separately.
\subsection{Expected Cholesky parameters}
Corollary~\ref{corr: tdisn} gives
\(\nu_i=\alpha_i-pa_i-2\) degrees of freedom for the marginal
\(t\) distribution of \(L_{\prec i]}\).  Consequently, if
\(\alpha_i>pa_i+3\),
\[
E(L_{\prec i]})=-U_{\prec i\succ}^{-1}U_{\prec i]}.
\]
If the stronger condition \(\alpha_i>pa_i+4\) holds, then
\begin{align}
\operatorname{Var}(L_{\prec i]})
&=\frac{U_{ii\mid\prec i\succ}}
{\alpha_i-pa_i-4}\,U_{\prec i\succ}^{-1},\label{eq:var-L}\\
E(D_{ii})
&=\frac{U_{ii\mid\prec i\succ}}
{\alpha_i-pa_i-4}.\label{eq:mean-D}
\end{align}
The Cholesky columns are mutually independent, so their joint
covariance is block diagonal with the blocks in \eqref{eq:var-L}.
The diagonal of \(L\) is deterministically one, and all non-parent
entries are deterministically zero.  These existence conditions are
stronger than the properness condition
\(\alpha_i>pa_i+2\).

\subsection{Mode in Cholesky coordinates}
We compute the joint density mode with respect to Lebesgue measure on
the free coordinates \((D_{ii},L_{\prec i]})\).  This base measure is
part of the definition: a nonlinear transformation of this point is
not, in general, a mode under the transformed coordinate measure.
From \eqref{D-parameterp}, the density is proportional to
 \[
 \prod_{i\in V}
 \lambda_i^{-\alpha_i/2}
 \exp\left\{-\frac{1}{2\lambda_i}
 \left[
 U_{ii\mid\prec i\succ}
 +(\ell_i+U^{-1}_{\prec i\succ}U_{\prec i]})^\top
 U_{\prec i\succ}
 (\ell_i+U^{-1}_{\prec i\succ}U_{\prec i]})
 \right]\right\}.
 \]
For each \(i\), first maximizing over \(\ell_i\) and then
differentiating with respect to \(\lambda_i\) gives the unique joint
mode
 \[
 \mathop{\times}_{i=1}^{p}
 \left(
 \frac{U_{ii\mid\prec i\succ}}{\alpha_i},
 -U^{-1}_{\prec i\succ}U_{\prec i]}
 \right).
 \]
The regression-coordinate component is the negative of the displayed
Cholesky component.

\begin{proposition}\label{prop:chol-mode}
Let \(Y_1,\ldots,Y_n\) be i.i.d. centered Gaussian observations and set
\(T=nS=\sum_rY_rY_r^\top\).  Under the prior
\(\pi^{\Xi_{\D}}_{U,\alpha}\), the posterior Cholesky-coordinate mode is
\[
 \mathop{\times}_{i=1}^{p}\left(
 \frac{(T+U)_{ii\mid\prec i\succ}}{\alpha_i+n},
 -(T+U)_{\prec i\succ}^{-1}(T+U)_{\prec i]}
 \right).
\]
Its regression coefficient is
\((T+U)_{\prec i\succ}^{-1}(T+U)_{\prec i]}\).  Completing these local
parameters gives a transformed Cholesky-mode plug-in, not a
Lebesgue-density MAP in covariance or precision coordinates.
\end{proposition}
%%%%%%%%%%%%%%%
\subsection*{The Jacobian of the mapping $\left(L,D\right)\mapsto\left(LD^{-1}L^{t}\right)^{E}$}
To derive the density of $\pi_{U,\alpha}^{\mathrm{R}_{\D}}$ we need to compute the Jacobian of the mapping 
\[
\psi\equiv \left(\left(L,D\right)\mapsto\left(LD^{-1}L^{t}\right)^{E}\right):\Theta_{\D}\rightarrow\mathrm{R}_{\D}.
\]
 The  Jacobian  of $\psi$ is a variant of similar transformations found in \cite{roverato00*, khare09a*}. For completeness we still compute this Jacobian  in the following lemma.

\begin{lemma} \label{lem:jacobian_of_psi_supp}
The Jacobian of the mapping $\psi:\left( (D, L\right)\mapsto\left( LD^{-1}L^t\right)^{E}$ is $\prod_{j=1}^p D_{jj} ^{-(pa_j+2)}$.
\end{lemma}

\begin{proof}

\noindent Let $\Upsilon \in \mathrm{R}_{\D}$, and $(D, L)\in \Theta_{\D}$  such that  $\widehat{\Upsilon}=LD^{-1}L^{\top}$. Note that for each $(i,j)\in E$,

\begin{equation} \label{Cholexpansion}
\Upsilon_{ij} = (LD^{-1}L^{\top})_{ij} = \sum_{k=1}^{p} L_{ik} L_{jk} D_{kk}^{-1} = \sum_{k=1}^j L_{ik} L_{jk} D_{kk}^{-1},
\end{equation}

\noindent
since $L$ is lower triangular. Now from  \eqref{Cholexpansion} it follows by noting that $L_{jj} = 1, \forall j $,
$$
\frac{\partial}{\partial L_{ij}} (LD^{-1}L^{\top})_{ij} = D^{-1}_{jj}, \;\;(i,j) \in E, \quad \quad \frac{\partial}{\partial D_{ii}} (LD^{-1}L^{\top})_{ii} = -D^{-2}_{ii}, \;\;i=1,2,\cdots,p.
$$

\noindent
Arrange the entries of $(D, L)\in { \Theta}_{\D}$ as $D_{11}$, $\{L_{2k}: \;$$
(2,k) \in E, 1 \leq k < 2\}$, $D_{22}$, $\{L_{3k}: \; (3,k) \in E$, $1 \leq k < 3\}$, $\ldots$,
$D_{p-1,p-1}$, $\left\{L_{pk}: \; (p,k) \in E, 1 \leq k < p\right\}, D_{pp}$,
and the entries of $\Upsilon \in \mathrm{R}_{\D}$ as $\Upsilon_{11}$, $\left\{\Upsilon_{2k}: \; (2,k)
\in E, 1 \leq k < 2\right\}$, $\Upsilon_{22}$, $\left\{\Upsilon_{3k}: \; (3,k) \in E, 1 \leq k <
3\right\}$, $\ldots$, $\Upsilon_{p-1,p-1}$, $\left\{\Upsilon_{pk}: \; (p,k) \in E, 1 \leq k < p\right\}$, $\Upsilon_{pp}$.
From \eqref{Cholexpansion} it is easily seen that $\Upsilon_{ij}$ depends on
 \[\left\{L_{jk}: \; (j,k) \in E, 1 \leq k < j\right\},\; \left\{ L_{ik}: \; (i,k) \in E, 1 \leq k < j \right\} \: \text{and}
  \: \left\{D_{kk}, 1 \leq k \leq j\right\}.
 \]
  \noindent Hence it is clear that $\Upsilon_{ij}$
is functionally independent of elements of $\Theta_{\D}$ that follow it in the arrangement described
above. Hence the gradient matrix of $\psi$ (with this arrangement) is a lower triangular matrix,
and the Jacobian of $\psi$ is therefore given as
$$
\prod_{i=1}^{p} \left(\prod_{j\in ch(i)} D_{jj}^{-1}\right)\prod^{p}_{i=1}D^{-2}_{ii}.
$$
\noindent
It follows from the expression above that the Jacobian of  $\psi$ is
\[
\prod_{j=1}^{p} D_{jj}^{-(pa_j+2)}.
\]
\end{proof}
%

%%%%%%%%%%%%%%%%%%%%%%%%%%%%%%%%%%%%%%%%%%%%%%%
%\section{Supplemental Section C}
%
%%---------------------------------------------------------- 

\subsection*{Exponential-family representation}\label{sub:curved}
Fix \(\alpha_i>pa_i+2\).  We consider the identifiable family
\[
\left\{\pi_{U,\alpha}^{R_{\D}}:U\in\mathrm{PD}_{\D}\right\},
\]
equivalently parameterized by \(U^E\in S_{\D}\) and the unique DAG
completion \(U\).  Thus, if
\(\pi_{U_1,\alpha}^{R_{\D}}=\pi_{U_2,\alpha}^{R_{\D}}\) almost
everywhere, then \(U_1^E=U_2^E\).  The ambient exponential-family
representation explains the distinction between perfect and
non-perfect DAGs.

\begin{lemma}\label{lem:exp_RG}
 Let $\D$ be DAG and let $\alpha$ be given. If $\D$ is perfect, then the Wishart family
\[ \left\{\pi^{\mathrm{R}_{\D}}_{U^{E},\alpha}: U^{E}\in \mathrm{S}_{\D}\right\}, \: \text{or equivalently} \:   \left\{\pi^{\mathrm{P}_{\D}}_{U^{E},\alpha}: U^{E}\in \mathrm{S}_{\D}\right\}
\]
 is a full regular exponential family in its natural affine support.
 If $\D$ is not perfect, then
 \(\{\pi^{\mathrm{R}_{\D}}_{U^{E},\alpha}:U^{E}\in
 \mathrm{S}_{\D}\}\) has a curved exponential-family representation in
 the larger moral-graph natural-parameter space.
 \end{lemma}
 \begin{proof}
If \(\D\) is perfect, its moral graph adds no edges and every completed
precision matrix has zeros outside \(E\).  Hence
\[
\tr(\widehat\Upsilon U)
=\langle\Upsilon^0,(U^E)^0\rangle_{Z_{\D}},
\]
and \(U^E\) ranges over an open subset of \(Z_{\D}\).  This gives a
full regular exponential family.

For a non-perfect DAG, moralization adds at least one edge.  The same
density can be embedded in the exponential family on
\(Z_{\D^{\mathrm m}}\) with canonical statistic
\((\widehat\Upsilon)^{E^{\mathrm m}}\) and natural parameter
\(-\tfrac12U^{E^{\mathrm m}}\).  The constraint
\(U\in\mathrm{PD}_{\D}\) makes the added moral-edge coordinates smooth,
nonlinear functions of the \(|E|\) free coordinates \(U^E\).  The
resulting natural-parameter set is therefore a \(|E|\)-dimensional
curved submanifold of the higher-dimensional moral-graph parameter
space.
\end{proof}

For non-perfect \(\D\), this DAG-constrained family is a proper
subfamily of the ambient family obtained by allowing an arbitrary
positive-definite \(U\).
%%%%%%%%%%%%%%%%%
%%%%%%%%%%%%%%%
\subsection*{The inverse DAG-Wishart for homogeneous DAGs}\label{sub:wh}
%%%%
We next show that the inverse DAG--Wishart family
\(\pi_{U,\alpha}^{\mathrm{S}_\D}\) contains the inverse Wishart subclass
introduced by Khare and Rajaratnam \cite{khare09a*} for Gaussian
covariance graph models.  For the associated special DAGs, the density
simplifies considerably.  Recall that a Gaussian covariance graph
model over an undirected graph \(\G=(V,\mathscr E)\) is denoted by
\(\mathscr N(\G_{\mathrm{cov}})\) and is defined as follows.
\begin{definition}Let $\mathrm{PD}_{\G_{\mathrm cov}}$ denote the set of positive definite matrices $\Sigma\in \mathrm{PD}_{p}(\R)$ such that $\Sigma_{ij}=0$ whenever $i\not\sim_{\G}j$, i.e., when $i$ and $j$ are not neighbors. Then the Gaussian covariance graph model over $\G$  is defined by $\mathscr{N}(\G_{\mathrm{cov}})=\left\{ \mathrm{N}_{p}(0,\Sigma): \Sigma\in \mathrm{PD}_{\G_{\mathrm{cov}}}\right\}$.
\end{definition}
A formal comparison between the DAG Wishart priors introduced in this paper and the covariance Wishart priors introduced in \cite{khare09a*} requires a few technical definitions.
\begin{definition}\label{def:homog_property}
\begin{itemize}
\renewcommand{\labelitemi}{$a)$}
\item A DAG $\D$ is called a homogeneous DAG of type I if it is transitive (i.e., $i\rightarrow j\rightarrow k$ implies that $i\rightarrow k$), and perfect. A DAG $\D$ is called a homogeneous DAG of type II if it is transitive and does not contain any induced subgraph of the form $j\leftarrow i\rightarrow k$.
\renewcommand{\labelitemi}{$b)$}
\item An undirected graph $\G=(V, \mathscr{E})$ is called homogeneous if 
\smallskip
$i\sim_{\G}j \Longrightarrow \mathrm{ne}(i)\cup \left\{ i\right\}\subseteq\mathrm{ne}(j)\cup \left\{j\right\}\:\: \text{or}\:\: \mathrm{ne}(j)\cup \left\{j \right\}\subseteq\mathrm{ne}(i)\cup\left\{ i\right\}$, for every $i,j\in V$.
\end{itemize}
\noindent Equivalently, a graph $\G$ is said to be homogeneous if it is decomposable and does not contain the $A_4$ path as an induced subgraph. The reader is referred to \cite{letac07*} for further details on homogeneous graphs.
\end{definition}
Note that if $\D$ is a homogeneous DAG of either types, then $\D^{u}$ is homogeneous.
On the other hand, if $\G=(V, \mathscr{E})$ is homogeneous, then one can construct a homogeneous DAG of type I or II that is a DAG version of $\G$. This can be achieved by using the Hasse tree associated with the homogeneous (undirected) graph and using the given orientation to obtain a DAG of type I. Reversing the orientation (i.e., redirecting all the arrows to the root of the tree) will yield a DAG of type II. More precisely we shall now show an example that constructs a DAG version that is homogeneous of type II. Let $\D$ be a directed version of $\G$ obtained by directing each edge $i\sim_{\G} j$ to a directed edge  $i\rightarrow j$ if $\mathrm{ne}(i)\cup \left\{i \right\}\subsetneq \mathrm{ne}(j)\cup \left\{ j\right\}$, or $j\rightarrow i$ if $\mathrm{ne}(j)\cup\left\{ j\right\}\subsetneq \mathrm{ne}(i)\cup \left\{ i\right\}$. If $\mathrm{ne}(i)\cup \left\{i\right\}=\mathrm{ne}(j)\cup \left\{ j\right\}$, an arbitrary direction is chosen. From Definition \ref{def:homog_property} one can check that $\D$ 
is a transitive DAG and it does not contain any
induced subgraph of the form $j\leftarrow i\rightarrow k$. In general, it can be shown that if $\D$ is a homogeneous DAG of type II and a DAG version of $\G$, then $\mathscr{N}({\D})$ is identical to the Gaussian covariance model $\mathscr{N}(\G_{\mathrm{cov}})$  in the sense that $\mathrm{PD}_{\G_{\mathrm{cov}}}=\mathrm{PD}_\D$ (see \cite{pearl93*} for instance for more details.) It is also evident, from the Markov equivalence of perfect DAGs and decomposable graphs, that for a homogeneous DAG $\D$ of type I which is a DAG version of $\G$, we have $\mathrm{PD}_{\G}=\mathrm{PD}_\D.$
 \begin{proposition}\label{prop:density_SG_hom}
 Let $\D=(V, E)$ be a homogeneous DAG of either type I or II and let $\G=(V, \mathscr{E})$ be a homogeneous graph.
 \begin{itemize}
 \renewcommand{\labelitemi}{$a)$}
 \item The density of $\pi_{U,\alpha}^{\mathrm{S}_\D}$  is given by \\
 \quad \;\; $z_{\D}(U,\alpha)^{-1}\exp\left\{-\frac{1}{2}\tr(\Sigma(\Gamma)^{-1} U)\right\}\prod_{i=1}^{p} \Sigma_{ii|\prec i\succ}^ {-\frac{1}{2}(\alpha_i+2ch_{i}(\D))}$, \quad \;\; where $ch_{i}(\D)=|\mathrm{ch}_\D(i)|$.
  \renewcommand{\labelitemi}{$b)$}
 \item If $\D$ is of type II and  a DAG version of $\G$, then the open cone $\mathrm{PD}_{\G_{\mathrm{cov}}}$ can be identified with $\mathrm{S}_\D$ via the bijective mapping
\begin{equation}\label{eq:2}
\left(\Gamma\mapsto\left(\Gamma\right)^{0}=\Sigma\left(\Gamma\right)\right):\mathrm{S}_\D\rightarrow \mathrm{PD}_{\G_{\mathrm{cov}}}.
\end{equation}
\noindent Let $\pi_{U,\alpha}^{\mathrm{PD}_{\G_{\mathrm{cov}}}}$ denote the probability image of the inverse DAG Wishart $\pi_{U,\alpha}^{\mathrm{S}_\D}$ under the mapping in  \eqref{eq:2}. Then the density of  $\pi_{U,\alpha}^{\mathrm{PD}_{\G_{\mathrm{cov}}}}$ w.r.t. Lebesgue measure is given by the expression in part(a) above.
 \end{itemize}
 \end{proposition}
 \begin{proof}~$a)$ It suffices to prove that for every $\Sigma\in\mathrm{PD}_\D$,
 \begin{equation}\label{eq:essential_eq}
\prod_{i\in V}\det(\Sigma_{\prec i\succ})= \prod_{i\in V}\Sigma_{ii| \prec i\succ}^{ch_{i}(\D)}.
\end{equation}
\begin{itemize}
\renewcommand{\labelitemi}{1)}
\item Suppose that $\D$ is homogeneous of type I. We shall first show that for every $i\in V$
\begin{equation}\label{eq:1}
\det(\Sigma_{\prec i\succ})=\prod_{\ell\in \mathrm{pa}(i)}\Sigma_{\ell\ell|\prec \ell\succ}.
\end{equation}
If $\mathrm{pa}(i)=\emptyset$ for some $i$, then by our convention $\det(\Sigma_{\prec i\succ})=1$  and $\Sigma_{\ell\ell|\prec \ell\succ}=1$ for any $\ell\in\mathrm{pa}(i)$ and therefore  \eqref{eq:1} holds. Now let $\ell_{0}$ be the smallest integer in $\mathrm{pa}(i)$. One then can easily check that since $\D$ is both transitive and perfect we have $\mathrm{pa}(i)=\left\{\ell_{0}\right\}\cup\mathrm{pa}(\ell_{0})$. From this we write $\det(\Sigma_{\prec i\succ})=\Sigma_{\ell_{0}\ell_{0}|\prec \ell_{0}\succ}\det(\Sigma_{\prec\ell_{0}\succ})$. Now by repeating this procedure we obtain the result in \eqref{eq:1}. Finally we write 
$$\prod_{i\in V}\det(\Sigma_{\prec i\succ})=\prod_{i\in V}\prod_{\ell\in\mathrm{pa}(i)}\Sigma_{\ell\ell|\prec \ell\succ}=\prod_{i\in V}\Sigma_{ii|\prec i\succ}^{ch_{i}(\D)}.$$
\renewcommand{\labelitemi}{2)}
\item Suppose $\D$ is homogeneous of type II. We shall proceed by induction. It is clear that  \eqref{eq:essential_eq} holds when $p=|V|=1$. Now by the inductive hypothesis assume that  \eqref{eq:essential_eq} holds for every homogeneous DAG of type II, connected or disconnected, with fewer vertices than $p=|V|$. Using the inductive hypothesis we shall show that  \eqref{eq:essential_eq} will also hold for $\D$ with $p$ vertices. Now let $\Sigma\in \mathrm{PD}_\D$ be given.
\end{itemize}
\begin{itemize}
\renewcommand{\labelitemi}{Case 1)}
\item Suppose that $\D$ is connected. Let $\D_{[1]}$ be the induced DAG on $V\setminus\left\{1\right\}$. It is clear that $\D_{[1]}$ is a homogeneous DAG of type II and therefore by the induction hypothesis $\prod_{i=2}^{p}\det(\Psi_{\prec i\succ})= \prod_{i=2}^{p}\Psi_{ii| \prec i\succ}^{ch_{i}(\D_{[1]})}$, where $\Psi=\Sigma_{V\setminus\left\{1\right\}}$. Note that $\D_{[1]}$ is an ancestral subgraph of $\D$ and hence $\mathrm{fa}_{{\mathcal D}_{[1]}}(i)=\mathrm{fa}_\D(i)$ for each $i=2,\ldots, p$ and consequently $\Psi_{\prec i \succ}=\Sigma_{\prec i \succ}$ and $\Psi_{ii| \prec i\succ}=\Sigma_{ii|\prec i\succ}$. All together these imply that
$\prod_{i=2}^{p}\det(\Sigma_{\prec i\succ})= \prod_{i=2}^{p}\Sigma_{ii| \prec i\succ}^{ch_{i}(\D_{[1]})}$. Now we claim that $\mathrm{fa}_\D(1)=V$. Assume to the contrary that $V\setminus \mathrm{fa}_\D(1)\neq \emptyset$. Since $\D$ is connected, this implies that there exist vertices $i\in \mathrm{fa}_\D(1)$  and $j\in V\setminus \mathrm{fa}_\D(1)$ such that $i,j$ are adjacent in $\D$. But this implies $j\rightarrow i\rightarrow 1$ or $j\leftarrow i\rightarrow 1$. By definition these induced subgraphs cannot occur in $\D$. Thus $\preceq 1\succeq =V$ and therefore we have $\det(\Sigma_{\prec 1\succ})=\Sigma_{11|\prec 1\succ}^{-1}\det(\Sigma)=\prod_{i=2}^{p}\Sigma_{ii|\prec i\succ}$.
Also the fact that $\mathrm{fa}_\D(1)=V$  implies that for each $i\in V\setminus\left\{1\right\}$ we have $\mathrm{ch}_{i}(\D_{[1]})=\mathrm{ch}_{i}(\D)-1$. Therefore
\begin{equation*}
\prod_{i\in V}\det(\Sigma_{\prec i\succ})=\det(\Sigma_{\prec 1\succ})\prod_{i=2}^{p}\det(\Sigma_{\prec i\succ})
=\prod_{i=2}^{p}\Sigma_{ii|\prec i\succ}\prod_{i=2}^{p}\Sigma_{ii|\prec i\succ}^{ch_{i}(\D_{[1]})}
=\prod_{i\in V}\Sigma_{ii| \prec i\succ}^{ch_{i}(\D)}.
\end{equation*}
%%%%%%%%%%%%%%%%
\renewcommand{\labelitemi}{Case 2)}
\item Suppose $\D$ is disconnected.  Let $\D_{1}$  and $\D_{2}$ denote respectively the induced subgraphs of $\D$ on $\mathrm{fa}_{D}(1)$ and $V\setminus \mathrm{fa}_\D(1)$. It is clear that $\D_{1}$ and $\D_{2}$  are both homogeneous of type II. In addition it is also easily verified that they are ancestral. Now let $\Psi=\Sigma_{\preceq 1\succeq}\in \mathrm{PD}_{\D_{1}}$ and 
$\Psi'=\Sigma_{V\setminus \mathrm{fa}_\D(1)}\in \mathrm{PD}_{\D_{2}}$. Now applying the induction hypothesis and the fact that $\D_{1}$  and $\D_{2}$ are disjoint we have:
\begin{align*}
\prod_{i\in V}\det(\Sigma_{\prec i\succ})&=\prod_{i\in \mathrm{fa}_\D(1)}\det(\Sigma_{\prec i\succ})  \prod_{i\in V\setminus\mathrm{fa}_\D(1)}\det(\Sigma_{\prec i\succ})\\
&=\prod_{i\in \mathrm{fa}_\D(1)}\det(\Psi_{\prec i\succ})  \prod_{i\in V\setminus\mathrm{fa}_\D(1)}\det(\Psi_{\prec i\succ}')\\
&=\prod_{i\in \mathrm{fa}_\D(1)}\det(\Psi_{ii|\prec i\succ})^{ch_{i}(\D_{1})}  \prod_{i\in V\setminus\mathrm{fa}_\D(1)}\det(\Psi_{ii|\prec i\succ}')^{ch_{i}(\D_{2})}\\
&=\prod_{i\in \mathrm{fa}_\D(1)}\det(\Sigma_{ii|\prec i\succ})^{ch_{i}(\D)}  \prod_{i\in V\setminus\mathrm{fa}_\D(1)}\det(\Sigma_{ii|\prec i\succ}')^{ch_{i}(\D)}\\
&=\prod_{i\in V}\det(\Sigma_{ii|\prec i\succ})^{ch_{i}(\D)}.
\end{align*}
\end{itemize}
$b)$  It is clear that the mapping in  \eqref{eq:2} is a diffeomorphism  and the Jacobian of this mapping is $1$. Thus  the functional form of the density $\pi_{U,\alpha}^{\mathrm{PD}_{\G_{\mathrm{cov}}}}$ w.r.t. Lebesgue measure is same as $\pi_{U,\alpha}^{\mathrm{S}_{\D}}$ given by Proposition \ref{prop:density_SG_hom}.
 \end{proof}
 \begin{Rem}
For a homogeneous graph $\G$, the distribution
$\pi_{U,\alpha}^{\mathrm{PD}_{\mathcal G_{\mathrm{cov}}}}$ with the density
in Proposition~\ref{prop:density_SG_hom} coincides with the covariance
Wishart family introduced by Khare and Rajaratnam \cite{khare09a*}.
\end{Rem}

%\bibliography{mybib}{}
%\bibliographystyle{natbib}

%%%%%%%%%%%%%%%%
\section*{Supplemental Section C: The DAG-Wishart on $\rm{P}_{\D}$ and its density w.r.t. Hausdorff measure}

\subsection*{Introduction}
In this section we define the DAG--Wishart distribution directly on the
precision-matrix model $\mathrm{P}_{\D}$ for an arbitrary DAG $\D$.  At the
level of probability measures, $\pi_{U,\alpha}^{\mathrm{P}_{\D}}$ is simply
the pushforward of $\pi_{U,\alpha}^{\Theta_{\D}}$ under
$(D,L)\mapsto LD^{-1}L^{\top}$.  If $\D$ is perfect, then
$\mathrm{P}_{\D}$ is an open subset of
$\mathrm{Z}_{\D}\cong\mathbb{R}^{|E|}$, so this pushforward has an ordinary
Lebesgue density.  For a non-perfect DAG, however,
$\mathrm{P}_{\D}$ is an $|E|$-dimensional embedded manifold in the larger
space $\mathrm{Z}_{\D^{\mathrm m}}$ and therefore has ambient Lebesgue
measure zero.  The appropriate reference measure is instead the
$|E|$-dimensional Hausdorff measure induced by a specified Euclidean metric.

 %------------------------------------------------------------
 \subsection*{Lebesgue measure of $\mathrm{P}_\D$}

Lemma~\ref{dfimoralf} gives
$\mathrm{P}_{\D}\subset\mathrm{P}_{\D^{\mathrm m}}
\subset\mathrm{Z}_{\D^{\mathrm m}}$.  For a non-perfect DAG, the model has
Lebesgue measure zero in every ambient linear subspace of
$\mathrm S_p(\mathbb R)$ that contains it.  The next lemma makes the
dimension argument precise.
 
\begin{lemma}\label{lem:not_open}
Suppose $\D=(V,E)$ is a non-perfect DAG and $\mathbb{E}$ is a linear
subspace of $\mathrm{S}_p(\mathbb{R})$ containing $\mathrm{P}_{\D}$.  Then
$\mathbb{E}$ contains $\mathrm{Z}_{\D^{\mathrm m}}$.  Consequently,
$\mathrm{P}_{\D}$ has Lebesgue measure zero in $\mathbb{E}$.
  \end{lemma}
  \begin{proof}
 For each $(i, j)\in E^{\mathrm{m}}$  with $j\leq i$ let us define the elementary symmetric matrix  $\widetilde{E}^{(ij)}\in \mathrm{S}_{p}(\R)$ as follows:
  \[
   \widetilde{E}^{(ij)}_{uv}=
  \begin{cases}
 1& \text{if  $\{u, v\}=\{i,j\}$, }\\
 0& \text{otherwise.}
  \end{cases}
  \]
The matrices $\widetilde{E}^{(ij)}$ form a basis of
$\mathrm{Z}_{\D^{\mathrm m}}$.  First, $\mathbb{E}$ contains
$\mathrm{Z}_{\D}$.  Indeed, positive diagonal precision matrices give the
diagonal directions by taking linear combinations, and, for each directed
edge $i\to j$, matrices of the form
$(I+t e_i e_j^{\top})(I+t e_i e_j^{\top})^{\top}$, together with the already
obtained diagonal directions, give $\widetilde E^{(ij)}$.  It remains to
obtain the moral edges.  Let $(i,j)\in E^{\mathrm m}\setminus E$, with
$i>j$.  Then some $k<j<i$ satisfies $i\to k\leftarrow j$.  Define the lower
triangular matrix $L^{(ij)}\in\mathcal{L}_{\D}$ by
  \[
      L^{(ij)}_{uv}=
  \begin{cases}
 1& \text{if  $(u, v)=(i,k)$, }\\
 1&  \text{if  $(u, v)=(j,k)$, }\\
 1& \text{if  $u=v$, }\\
 0& \text{otherwise.}
  \end{cases} 
  \]
Then
$L^{(ij)}(L^{(ij)})^{\top}=T+\widetilde{E}^{(ij)}$ for some
$T\in\mathrm{Z}_{\D}$.  Hence $\widetilde{E}^{(ij)}\in\mathbb{E}$ and
$\mathrm{Z}_{\D^{\mathrm m}}\subseteq\mathbb{E}$.
 
The model $\mathrm{P}_{\D}$ is a smooth manifold of dimension $|E|$,
whereas $\dim(\mathrm{Z}_{\D^{\mathrm m}})=|E^{\mathrm m}|>|E|$ for a
non-perfect DAG.  Thus $\dim(\mathbb{E})>|E|$, and the embedded
$|E|$-dimensional manifold $\mathrm{P}_{\D}$ has Lebesgue measure zero in
$\mathbb{E}$.
\end{proof}
Consequently, Lemma \ref{lem:not_open}  implies that  if $\D$ is non-perfect then $\pi^{\mathrm{P}_{\D}}_{U,\alpha}$  has no density w.r.t. Lebesgue measure. 
%---------------------------------------------------------- 
\subsection*{ The density of  $\pi_{U,\alpha}^{\mathrm{P}_{\D}}$ w.r.t. Hausdorff measure} 

We now derive the density of $\pi_{U,\alpha}^{\mathrm{P}_{\D}}$ with
respect to Hausdorff measure.  Write $q=|E|$, where the diagonal pairs are
included in $E$, and let
\[
 \mathcal O_{\D}=(0,\infty)^p\times\mathbb{R}^{q-p}.
\]
For $\theta=(d_1,\ldots,d_p,(\ell_{ij})_{(i,j)\in E,\,i>j})\in
\mathcal O_{\D}$, set
\[
 D(\theta)=\operatorname{diag}(d_1,\ldots,d_p),\qquad
 L(\theta)=I+\sum_{(i,j)\in E,\,i>j}\ell_{ij}e_i e_j^{\top}.
\]
Thus $\Theta_{\D}$ is identified with the open subset
$\mathcal O_{\D}\subset\mathbb{R}^{q}$; it is not itself a vector space
because $D$ must be positive and $L$ must have unit diagonal.  In the target,
we identify $\mathrm{Z}_{\D^{\mathrm m}}$ with
$\mathbb{R}^{|E^{\mathrm m}|}$ by retaining one lower-triangular coordinate
for each symmetric position.  All Euclidean norms and Hausdorff measures
below refer to these specified coordinates.  Using the Frobenius metric on
symmetric matrices would instead weight off-diagonal directions twice in
the Gram matrix, so the Jacobian and displayed density would have to be
adjusted.

The map
\[
 \psi:\mathcal O_{\D}\longrightarrow\mathrm{Z}_{\D^{\mathrm m}},
 \qquad \psi(\theta)=L(\theta)D(\theta)^{-1}L(\theta)^{\top},
\]
is smooth and injective.  Order the domain coordinates as
$D_{11},L_{21}$ (when $(2,1)\in E$), $D_{22},L_{31}$ (when
$(3,1)\in E$), and so on.  Let
$\mathscr I=E^{\mathrm m}\setminus E$.  In the target, order first the
positions in $E$ in the corresponding order and then the positions in
$\mathscr I$ lexicographically.  The required partial derivatives are
\vspace{-0.2cm}
\begin{equation}\label{derd}
\frac{\partial (LD^{-1}L^t)_{kl}}{\partial D_{ii}}= -D^{-2}_{ii}L_{ki}L_{li}
\end{equation}
\begin{equation}\label{derl}
\frac{\partial (LD^{-1}L^t)_{kl}}{\partial L_{ij}}=\delta_{ik}D_{jj}^{-1}L_{lj}+\delta_{il}D_{jj}^{-1}L_{kj},
\end{equation}
where $\delta_{uv}$ is the Kronecker delta.  Partition $D\psi$ into the
$q\times q$ block $A_{\psi}$ comprising the rows indexed by $E$ and the
$|\mathscr I|\times q$ block $C_{\psi}$ comprising the remaining rows.
The block $A_{\psi}$ is the Jacobian from Lemma
\ref{lem:jacobian_of_psi}; in particular it is nonsingular.  Hence $\psi$ is
an immersion, and the area formula gives the $q$-dimensional Jacobian
\begin{align*}
\mathbf{J}\psi(D,L)
&=\sqrt{\det\!\left((D\psi)^{\top}D\psi\right)}\\
&=\sqrt{\det(A_{\psi}^{\top}A_{\psi}+C_{\psi}^{\top}C_{\psi})}\\
&=|\det(A_{\psi})|\sqrt{\det(I+A_{\psi}^{-\top}C_{\psi}^{\top}C_{\psi}A_{\psi}^{-1})}\\
&=\prod_{j=1}^{p} D_{jj}^{-(pa_j+2)}
\sqrt{\det(I+A_{\psi}^{-\top}C_{\psi}^{\top}C_{\psi}A_{\psi}^{-1})}.
\end{align*}
Therefore we have proved the following.
\begin{theorem}\label{thm:hausdorff-density}
Let $A_{\psi}$ and $C_{\psi}$ be the derivative blocks defined above.
With respect to the restriction to $\mathrm P_{\D}$ of the
$q$-dimensional Hausdorff measure induced by the chosen lower-triangular
target coordinates, the density of
$\pi_{U,\alpha}^{\mathrm{P}_{\D}}$ is
\begin{equation}\label{hd}
z_{\D}(U,\alpha)^{-1}\exp\{-\tfrac{1}{2}\tr(\Omega U)\}
\prod_{i=1}^{p}D_{ii}^{-\alpha_i/2+pa_i+2}
\det(I+A_{\psi}^{-\top}C_{\psi}^{\top}C_{\psi}A_{\psi}^{-1})^{-1/2},
\end{equation}
where $D$ and $L$ are the unique modified-Cholesky coordinates of $\Omega$.
\end{theorem}
%%%%%
\begin{figure}
  \centering
  \ColliderDAG
  \caption{The non-perfect collider DAG $3\to1\leftarrow2$ used in the
  Hausdorff-density example.}
  \label{fig-33}
\end{figure}

%%%%%

\begin{Ex}
For the collider $3\to1\leftarrow2$ in Figure~\ref{fig-33}, use
\[
\begin{aligned}
\theta&=(D_{11},L_{21},D_{22},L_{31},D_{33}),\\
\eta&=(\Omega_{11},\Omega_{21},\Omega_{22},
       \Omega_{31},\Omega_{33},\Omega_{32}).
\end{aligned}
\]
Equations~\eqref{derd} and~\eqref{derl} give
\[
M_{\psi}=
\left(
\begin{matrix}
       -D_{11}^{-2}&  0&      0&         0&      0\\
	    -L_{21}D_{11}^{-2}&       D_{11}^{-1}&      0&         0&      0\\
  -L_{21}^2D_{11}^{-2}& 2L_{21}D_{11}^{-1}& -D_{22}^{-2}&         0&   0\\
  -L_{31}D_{11}^{-2}&         0&      0&       D_{11}^{-1} &     0\\
	   -L_{31}^2D_{11}^{-2}&         0&      0& 2L_{31}D_{11}^{-1}& -D_{33}^{-2}\\
-L_{21}L_{31}D_{11}^{-2}&   L_{31}D_{11}^{-1}&      0&   L_{21}D_{11}^{-1}& 0\\
\end{matrix}
\right)
\]
Consequently,
\[
\mathbf{J}\psi(D,L)=D_{11}^{-4}D_{22}^{-2}D_{33}^{-2}
\sqrt{(1+L_{21}^{2})(1+L_{31}^{2})}.
\]
Thus, under the stated coordinate metric, the density with respect to
$\mathcal{H}^{5}$ on $\mathbb{R}^{6}$ is
\[
z_{\D}(U,\alpha)^{-1}\exp\{-\tfrac{1}{2}\tr(\Omega U)\}
D_{11}^{-\alpha_1/2+4}D_{22}^{-\alpha_2/2+2}D_{33}^{-\alpha_3/2+2}
\{(1+L_{21}^{2})(1+L_{31}^{2})\}^{-1/2},
\]
where the $D_{ii}$ and $L_{ij}$ are the unique functions of $\Omega$
defined by the modified Cholesky factorization.
\end{Ex}

%\newpage

\section*{Supplemental Section D: Computational algorithms and more related results}

\subsection*{DAG-constrained MLE and Cholesky-mode plug-in}
\label{sec:suppAlgo}

\begin{algo}[Maximum Likelihood]\label{algo:mle}
Let $Y_1,\ldots,Y_n$ be independent centered Gaussian observations and
$S=n^{-1}\sum_{r=1}^nY_rY_r^{\top}$.  Assume
$n\geq\max_i(pa_i+1)$, which makes every displayed Schur complement
positive almost surely.  (If a mean is first estimated, the corresponding
sufficient condition is $n\geq\max_i(pa_i+2)$.)  For each $i$ set
\[
\lambda_i=S_{ii|\prec i\succ}\in\mathbb{R}_{+}\quad\text{and}\quad\boldsymbol{\beta}_i=S_{\prec i\succ}^{-1}S_{\prec i]}\in\mathbb{R}^{pa_i}.
\]
Note that $\lambda_i=S_{ii}$ whenever $\mathrm{pa}(i)=\emptyset$.  For
$i=p,p-1,\ldots,1$:
\begin{enumerate}
\item Initialize $\widehat{\Sigma}_{ii}=\lambda_i$ for each $i$ such that $\mathrm{pa}(i)=\emptyset$ (in particular for $i=p$);
\item set $\widehat{\Sigma}_{\prec i]}=\widehat{\Sigma}_{\prec i\succ}\boldsymbol{\beta}_i$ if $\mathrm{pa}(i)\ne \emptyset$;
\item set $\widehat{\Sigma}_{ii}=\lambda_i+\boldsymbol{\beta}_i^{\top}\widehat{\Sigma}_{\prec i\succ}\boldsymbol{\beta}_i$ if $\mathrm{pa}(i)\ne \emptyset$;
\item set $\widehat{\Sigma}_{\nprec i]}=\widehat{\Sigma}_{\nprec i\succ}\widehat{\Sigma}_{\prec i\succ}^{-1}\widehat{\Sigma}_{\prec i]}$ if $\mathrm{pa}(i)\ne \emptyset$, otherwise set $\widehat{\Sigma}_{\nprec i]}=0$. 
\end{enumerate}
\end{algo} 
The precision-matrix MLE is $\widehat\Omega=\widehat\Sigma^{-1}$.

Proposition~\ref{prop:chol-mode} gives the joint posterior mode with respect
to Lebesgue measure in the free modified-Cholesky coordinates.  Applying the
same covariance completion to its local residual variances and regression
coefficients gives the following plug-in.
\begin{algo}[Transformed Cholesky-mode plug-in]\label{algo:map}
Let $S=n^{-1}\sum_{r=1}^nY_rY_r^{\top}$.  For $i=1,\ldots,p$, set
\[
 \lambda_i=\frac{(nS+U)_{ii\mid\prec i\succ}}{\alpha_i+n},\qquad
 \boldsymbol\beta_i=(nS+U)_{\prec i\succ}^{-1}(nS+U)_{\prec i]}.
\]
When $\mathrm{pa}(i)=\emptyset$, use
$\lambda_i=(nS+U)_{ii}/(\alpha_i+n)$ and
$\boldsymbol\beta_i=0$.  For $i=p,p-1,\ldots,1$, apply steps 1--4 of
Algorithm~\ref{algo:mle}, writing the result as
$\widetilde\Sigma_{\mathrm{chol}}$.
\end{algo}
Set $\widetilde\Omega_{\mathrm{chol}}=
\widetilde\Sigma_{\mathrm{chol}}^{-1}$.  These matrices are transformations
of the Cholesky-coordinate mode.  They are not, in general, MAP estimators
with respect to Lebesgue measure in covariance or precision coordinates.

\subsection*{Historical covariance-estimation results}\label{sec:suppCov}

Table~\ref{tab:Estimation-p500-Covariance-detail} preserves the reported
covariance results for $p=500$ and admissible-edge probability $0.01$.
The notation is corrected to match Section~\ref{sec:empirical}:
$\bar\Sigma$ is the full posterior mean,
$\widehat\Omega_{\mathrm{coord}}^{-1}$ is the inverse of the completed
posterior-mean precision coordinates, and
$\widetilde\Sigma_{\mathrm{chol}}$ is the transformed Cholesky-mode plug-in.
Positive entries mean lower loss than the constrained MLE; negative entries
mean higher loss.  The different rankings under $L_1$ and $L_2$ illustrate
that none of the three procedures uniformly dominates.  Because the public
archive contains neither saved replicate-level output nor an aggregation
script, the values and boldface rankings are historical and have not been
independently reproduced; no Monte Carlo uncertainty is available.

\begin{table}[H]
\begin{center}
\begin{tabular}{ |c|c|cc|cc|cc| }
\hline
&&\multicolumn{2}{ |c| }{n=30}& \multicolumn{2}{ |c| }{n=50} &  \multicolumn{2}{ |c| }{n=100} \\
\hline
$(c,U)$&Estimator& $L_1$ & $L_2$ & $L_1$ & $L_2$ & $L_1$ & $L_2$\\ \hline
 \multirow{3}{*}{$(2.5,I(3))$} &$\bar\Sigma$             & -9.8\% & 4.6\% & -8.0\% & 1.4\%  & -5.1\% & 0.2\% \\
                                              & $\widehat\Omega_{\mathrm{coord}}^{-1}$  & 27.9\% & -113.2\%  & 17.4\% & -100.2\%  & 8.8\% & -68.7\%\\
                                              & $\widetilde\Sigma_{\mathrm{chol}}$             & 27.4\% & -32.0\%  & 17.4\% & -27.0\%  & 9.1\%   & -16.4\%\\
\hline
\multirow{3}{*}{$(3,I(3))$} &$\bar\Sigma$                & 1.0\% & 10.5\%  & -1.2\% & 5.4\%  & -1.4\% & 2.6\% \\
                                              & $\widehat\Omega_{\mathrm{coord}}^{-1}$  & 30.1\% & -130.6\%  & 19.2\% & -115.6\%  & 10.0\% & -79.4\%\\
                                              & $\widetilde\Sigma_{\mathrm{chol}}$             & 27.1\% & -45.2\%  & 17.0\% & -38.2\%  & 8.8\%   & -23.6\%\\
 \hline
 \multirow{3}{*}{$(3.5,I(3))$} &$\bar\Sigma$           & 7.6\% & {\bf12.3\%} & 4.0\% & {\bf 6.8\%}  & 1.6\% & {\bf 3.4\%} \\
                                              & $\widehat\Omega_{\mathrm{coord}}^{-1}$  & {\bf 31.0\%} & -148.1\%  & {\bf 19.9\%} & -131.7\%  & {\bf 10.5\%} & -90.8\%\\
                                              & $\widetilde\Sigma_{\mathrm{chol}}$             & 26.0\% & -58.8\%  & 15.9\% & -50.1\%  & 8.0\%   & -31.7\%\\
 \hline
  \multirow{3}{*}{$(3,I(2.5))$} &$\bar\Sigma$          &   7.9\% & 11.8\% & 4.5\% & 6.2\%  & 2.0\% & 2.9\% \\
                                              & $\widehat\Omega_{\mathrm{coord}}^{-1}$  & 30.8\% & -141.6\%  & 19.7\% & -124.3\%   & 10.4\% & -84.7\%\\
                                              & $\widetilde\Sigma_{\mathrm{chol}}$             & 24.6\% & -49.2\%  & 14.8\% & -42.1\%  & 7.3\%   & -26.6\%\\
 \hline
  \multirow{3}{*}{$(3,I(3.5))$} &$\bar\Sigma$          & -9.8\% & 8.2\% & -8.7\% & 4.0\%  & -5.9\% & 1.8\% \\
                                              & $\widehat\Omega_{\mathrm{coord}}^{-1}$  & 27.4\% & -120.3\%  & 16.9\% & -107.6\%  & 8.4\% & -74.5\%\\
                                              & $\widetilde\Sigma_{\mathrm{chol}}$             & 27.4\% & -41.6\%  & 17.4\% & -34.6\%  & 9.1\%   & -21.0\%\\
 \hline
 \end{tabular}
\end{center}
\caption{Historically reported relative changes from the constrained MLE
when estimating $\Sigma$, with $p=500$.  Positive values indicate lower
loss.  The three rows are the full posterior mean, the inverse of the
coordinate-completed precision estimate, and the transformed
Cholesky-mode plug-in.  Saved outputs and Monte Carlo uncertainty are
unavailable.}
\label{tab:Estimation-p500-Covariance-detail}
\end{table}

\subsection*{Sensitivity to edge density and to one contamination design}
\label{sec:suppOutlier}
\begin{figure}[H]
\begin{center}
\includegraphics[width=0.8\textwidth]{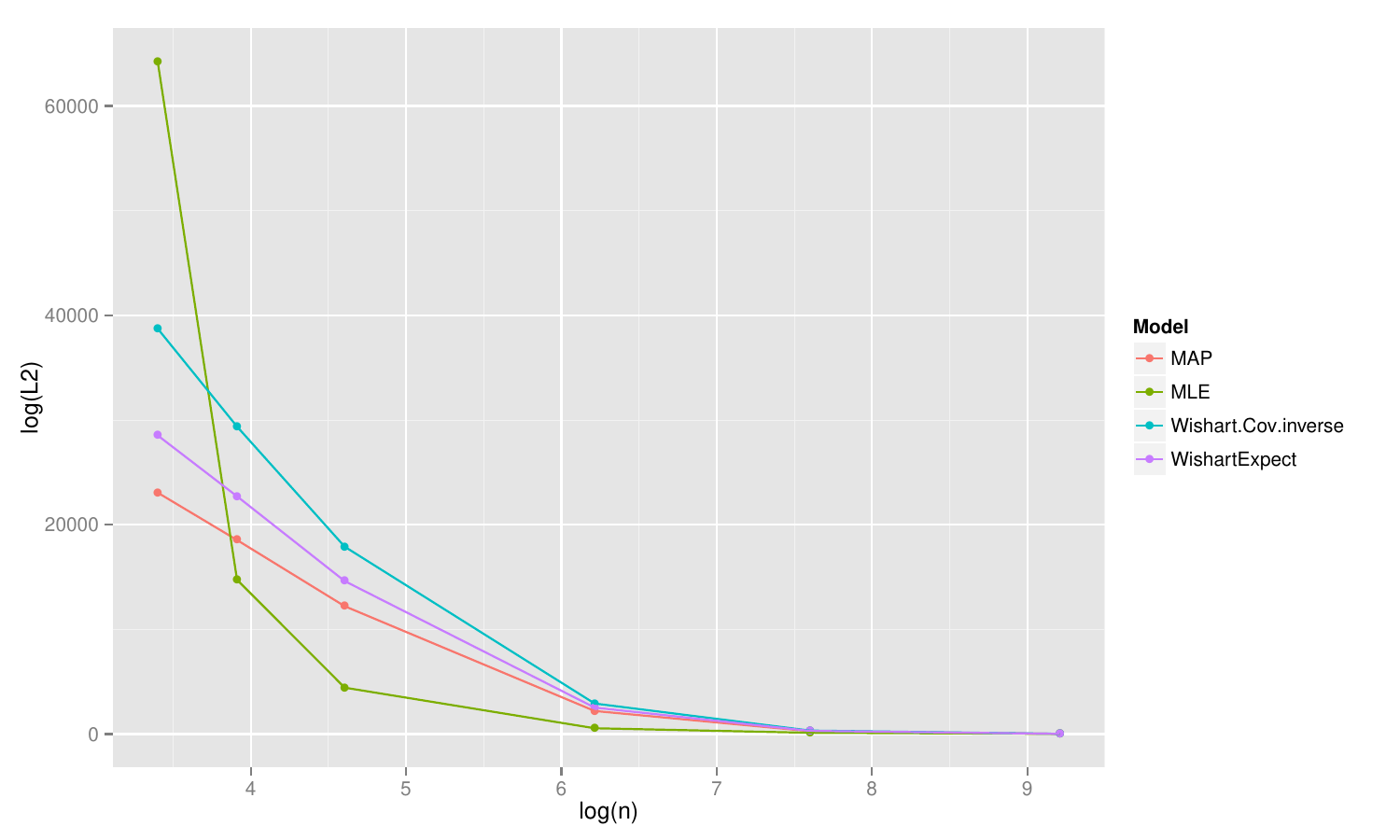}
\caption{Historical $L_2$-loss curves as $n$ increases, at admissible-edge
probability $0.015$.  The source values and uncertainty bands are
unavailable.}
\label{fig:nvserr}
 \end{center}
\end{figure}

Table~\ref{tab:Estimation-p500-Precision-Sparsity} preserves the historical
results obtained with $p=500$, $c=3$, $u=3$, and admissible-edge
probabilities from $0.005$ to $0.02$.  The results are strongly
design-dependent.  At edge probability $0.02$ and $n=100$, every listed
procedure is much worse than the constrained MLE under the reported $L_2$
loss, with relative changes between $-795.0\%$ and $-834.4\%$.  The table
therefore refutes any claim of uniform dominance and shows that the chosen
hyperparameters do not transfer reliably across edge densities.  It does
not identify an alternative setting, because that would require a separate,
properly replicated tuning study.  Figure~\ref{fig:nvserr} is the preserved
historical loss curve for edge probability $0.015$; its numerical source
data are absent from the archive.
 
\begin{table}[htbp]
\begin{center}
\small
\begin{tabular}{ |c|c|cc|cc|cc| }
\hline
&&\multicolumn{2}{ |c| }{n=30}& \multicolumn{2}{ |c| }{n=50} &  \multicolumn{2}{ |c| }{n=100} \\
\hline
Edge prob.&Estimator& $L_1$ & $L_2$ & $L_1$ & $L_2$ & $L_1$ & $L_2$\\ \hline
&&&&&&&\\
\multirow{3}{*}{0.005}        &$\widehat\Omega_{\mathrm{coord}}$          & 33.7\% & 72.7\%  & 21.4\% & 53.7\% & 11.3\% & 32.5\% \\
                                              & $\bar\Sigma^{-1}$& 38.9\% & 68.3\%  & 25.6\% & 51.0\% & 13.9\% & 31.8\%\\
                                              & $\widetilde\Omega_{\mathrm{chol}}$        & 22.0\% & 62.2\% & 11.8\% & 41.5\% & 5.1\% & 22.3\%\\
																							&&&&&&&\\

 \hline
&&&&&&&\\
\multirow{3}{*}{0.01}          &$\widehat\Omega_{\mathrm{coord}}$          & 39.2\% & 80.5\%  & 24.7\% & 60.5\% & 12.9\% & 34.6\% \\
                                              & $\bar\Sigma^{-1}$& 47.4\% & 65.9\%  & 31.1\% & 39.9\% & 16.7\% & 13.8\%\\
                                              & $\widetilde\Omega_{\mathrm{chol}}$        & 34.4\% & 81.5\% & 20.1\% & 62.3\% & 9.7\% & 37.7\%\\
																							&&&&&&&\\

 \hline
&&&&&&&\\
 \multirow{3}{*}{0.015}      &$\widehat\Omega_{\mathrm{coord}}$                & 45.1\% & 57.9\% & 23.2\% & -34.5\%  & 8.6\% & -244.4\% \\
                                              & $\bar\Sigma^{-1}$     & 49.1\% & 45.8\%  & 28.4\% & -67.0\%  & 8.3\% & -307.8\%\\
                                              & $\widetilde\Omega_{\mathrm{chol}}$             & 48.6\% & 64.3\%  & 26.9\% & -14.3\%  & 12.4\%   & -198.2\%\\
																							&&&&&&&\\
\hline
&&&&&&&\\
\multirow{3}{*}{0.02}    &$\widehat\Omega_{\mathrm{coord}}$                  & 38.1\% & 60.0\% & -7.0\% & -118.1\%  & -59.1\% & -812.6\% \\
                                              & $\bar\Sigma^{-1}$  & 36.3\% & 58.7\%  & -13.3\% & -124.7\%  & -68.4\% & -834.4\%\\
                                              & $\widetilde\Omega_{\mathrm{chol}}$          & 47.9\% & 60.8\%  & 6.4\% & -113.4\%  & -44.7\%   & -795.0\%\\
																							&&&&&&&\\
\hline

\end{tabular}
\end{center}
\caption{Historically reported relative changes from the constrained MLE
when estimating $\Omega$ at different admissible-edge probabilities, with
$p=500$ and $(c,u)=(3,3)$.  Positive values indicate lower loss and negative
values higher loss.  The procedures are
$\widehat\Omega_{\mathrm{coord}}$, $\bar\Sigma^{-1}$, and
$\widetilde\Omega_{\mathrm{chol}}$.  Saved outputs and Monte Carlo
uncertainty are unavailable.}
\label{tab:Estimation-p500-Precision-Sparsity}
\end{table}

The companion driver also calls \texttt{rmvDAG} with the arguments
\texttt{errDist="mixN100"} and \texttt{mix=0.01}, described in the original
study as replacing $1\%$ of Gaussian errors by $N(0,100)$ errors.  The
generated columns are then centered and standardized before estimation.
The exact semantics depend on the unrecorded package version.  The driver
requests 15 batches of 20 replicates for each sample size, but the resulting
\texttt{.Rda} file is absent.

Table~\ref{tab:RobustEstimation-p500-Precision} preserves the reported
relative changes for this single contamination mechanism.  Several
precision-loss changes are larger than in the uncontaminated table, while
some covariance $L_2$ changes remain negative.  Without the saved output,
Monte Carlo uncertainty, alternative contamination mechanisms, influence
analysis, or breakdown calculations, this experiment supports only a
design-specific sensitivity comparison.  It does not establish robustness
or a general guarantee under misspecification.

\begin{table}[H]
\begin{center}
\begin{tabular}{ |c|c|cc|cc|cc| }
\hline
&&\multicolumn{2}{ |c| }{n=30}& \multicolumn{2}{ |c| }{n=50} &  \multicolumn{2}{ |c| }{n=100} \\
\hline
Target&Estimator& $L_1$ & $L_2$ & $L_1$ & $L_2$ & $L_1$ & $L_2$\\ \hline
&&&&&&&\\
\multirow{3}{*}{$\Omega$} &$\widehat\Omega_{\mathrm{coord}}$& 65.1\% & 95.0\% & 50.9\% & 87.3\% & 32.2\% & 69.7\% \\
 & $\bar\Sigma^{-1}$& 72.1\% & 95.2\%  & 58.0\% & 89.7\% & 37.8\% & 75.0\%\\
 & $\widetilde\Omega_{\mathrm{chol}}$ & 59.7\% & 93.3\% & 44.9\% & 83.3\% & 27.0\% & 63.4\%\\
&&&&&&&\\
 \hline
&&&&&&&\\
\multirow{3}{*}{$\Sigma$}   &$\bar\Sigma$& 26.6\% & 9.2\% &  23.1\% & 4.6\% & 15.9\% & 1.6\% \\
                                                 & $\widehat\Omega_{\mathrm{coord}}^{-1}$& 45.1\% &-67.5\%  & 34.3\% & -40.2\% & 21.5\% &-16.4\%\\
                                                 & $\widetilde\Sigma_{\mathrm{chol}}$ & 17.7\% & -13.8\% & 28.6\% & -5.4\% & 17.7\% & -1.0\%\\% \hline
																								&&&&&&&\\
\hline
\end{tabular}
\end{center}
\caption{Historically reported relative changes from the constrained MLE in
the single \texttt{mixN100} contamination design, with $p=500$ and
$(c,u)=(3,3)$.  Positive values indicate lower loss.  Saved outputs,
package versions, and Monte Carlo uncertainty are unavailable; this is not
a general robustness assessment.}
\label{tab:RobustEstimation-p500-Precision}
\end{table}

\subsection*{Historical call-center prediction example}\label{sec:callcenter}
Following the descriptions in \cite{Levina08} and \cite{rajaratnam08*}, the
original analysis reports 239 operating days from a financial call center in
2002, with call counts $N_{ij}$ in 102 ten-minute intervals between 7 a.m.
and midnight.  It uses the variance-stabilizing transformation
$x_{ij}=\sqrt{N_{ij}+1/4}$.  The raw data file is not included in the public
archive, so these exclusions, counts, and transformations cannot be checked
in this revision.  Chronological order is imposed as the admissible parent
order.  This is a predictive modeling restriction, not evidence that the
estimated edges are causal or that latent common causes are absent.

The task is to predict the second 51 intervals from the first 51.  Write
$x_i=(x_i^{(1)},x_i^{(2)})$ and partition
\[
\mu=\begin{bmatrix}\mu^{(1)}\\ \mu^{(2)}\end{bmatrix},\qquad
\Sigma=\begin{bmatrix}
\Sigma_{11}&\Sigma_{12}\\
\Sigma_{21}&\Sigma_{22}
\end{bmatrix}.
\]
For an estimated mean and covariance, the conditional-mean predictor is
\[
\widehat x^{(2)}=\widehat\mu^{(2)}+
\widehat\Sigma_{21}\widehat\Sigma_{11}^{-1}
(x^{(1)}-\widehat\mu^{(1)}).
\]
The reported split uses the first 205 days for training and the last 34 for
testing.  The four procedures are the unrestricted sample-covariance
predictor (Naive-MLE), graph-constrained MLEs after LassoDAG and DAG--Wishart
selection, and the inverse
$\widehat\Omega_{\mathrm{coord}}^{-1}$ after DAG--Wishart selection.  The
last procedure is called ``DAG-W-Precision'' in the historical figure, but
it is not the posterior mean of $\Sigma$.

The script fixes $\kappa=0.1$ for LassoDAG and $c=1,b=3$ for DAG--Wishart;
these choices were not selected by cross-validation on the displayed test
split.  More importantly, after estimating and removing the training mean,
the script uses $n\operatorname{var}(x)$ and increments $\alpha$ by $n$.
Under the flat-mean analysis in Remark~\ref{rem:unknown-mean}, the update is
instead $(n-1)\operatorname{var}(x)$ and $\alpha+(n-1)\mathbf 1$.  Thus the
reported coordinate-shrinkage predictor does not implement the corrected
unknown-mean posterior.  The unavailable raw data and fitted-model object
prevent recomputation.

\begin{figure}[H]
\begin{center}
\includegraphics[width=0.4\textwidth]{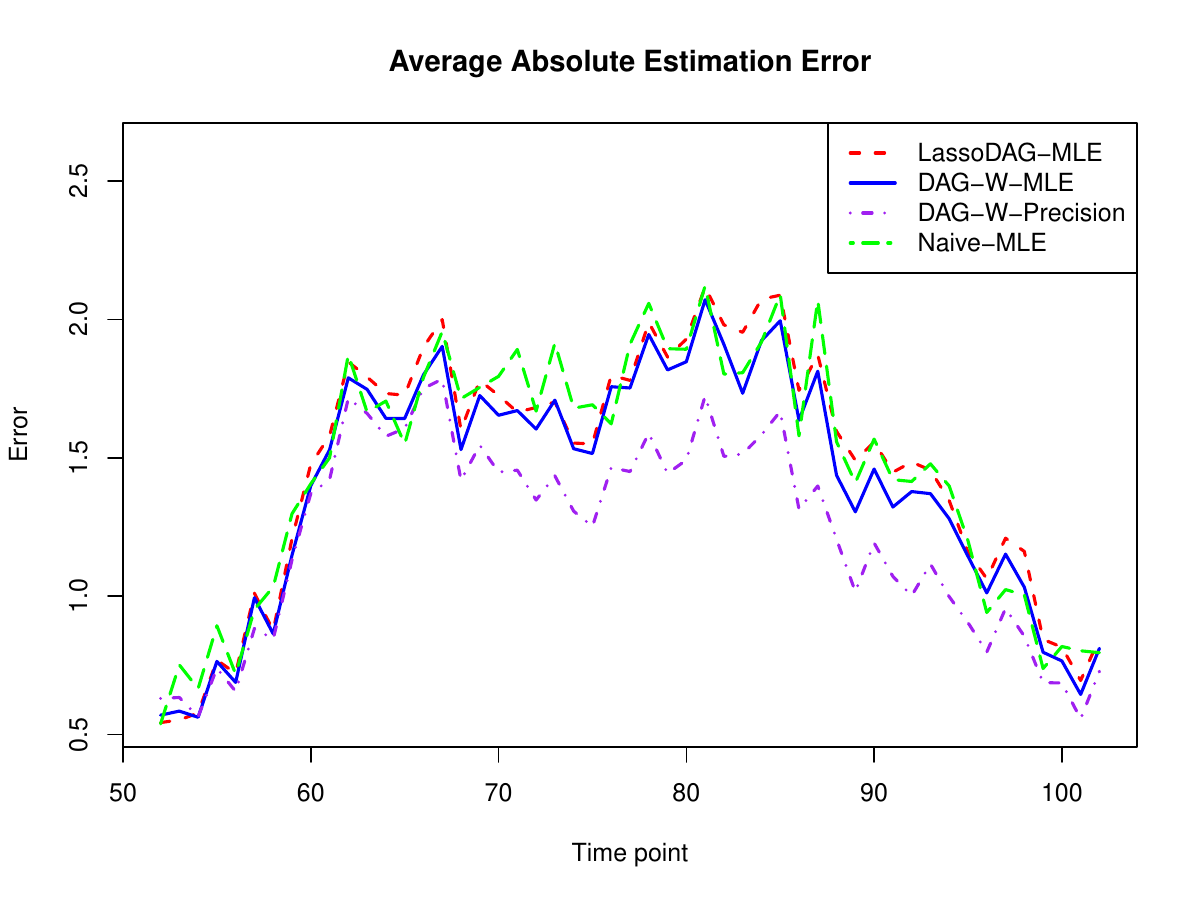}
\caption{Historical mean absolute errors.  ``DAG-W-Precision'' denotes
$\widehat\Omega_{\mathrm{coord}}^{-1}$; inputs and uncertainty are unavailable.}
\label{fig:CallCenter}
\end{center}
\end{figure}

For response interval $k=1,\ldots,51$, the plotted average absolute error is
\[
E_k=\frac1{34}\sum_{i=1}^{34}
\left|x_{i,51+k}-\widehat x_{i,51+k}\right|.
\]
Figure~\ref{fig:CallCenter} and Table~\ref{CallCenterMSE} preserve the
historically reported point estimates.  They have no resampling uncertainty
and cannot establish that one graph-selection method is generally superior.
The table statistic is the mean, over test days, of the 51-dimensional sum
of squared errors,
\[
\operatorname{meanSSE}=\frac1{34}\sum_{i=1}^{34}
\lVert x_i^{(2)}-\widehat x_i^{(2)}\rVert_2^2,
\]
not a per-coordinate mean squared error.
\enlargethispage{3\baselineskip}
\begin{table}[H]
\centering
\small
\begin{tabular}{|c|c|c|c|c|}
  \hline
 & Naive-MLE & LassoDAG-MLE & DAG-W-MLE &DAG-W-coordinate\\ 
  \hline
Mean SSE & 172.976 &  166.138 & 142.730 & 123.438 \\ 
   \hline
\end{tabular}
\caption{Historically reported mean test-day sums of squared prediction
errors for the call-center example.  ``DAG-W-coordinate'' uses
$\widehat\Omega_{\mathrm{coord}}^{-1}$.  These values are not independently
verified and do not include uncertainty.}
\label{CallCenterMSE}
\end{table}

%%%%%%%%%%%%%%%%%%%%%%%%%%%%%%%%%%%%%%%%%%%%%%%%%%%%%%%%%%%%%%

\end{document}